\documentclass[a4paper,12pt]{article}
\usepackage{amssymb}
\usepackage{amsthm}
\usepackage{amsmath}
\usepackage{latexsym}
\usepackage[dvips]{graphicx}
\usepackage{exscale}
\usepackage[latin1]{inputenc}

\newtheorem{theorem}{Theorem}[section]

\newtheorem{lemma}{Lemma}[section]
\newtheorem{proposition}{Proposition}[section]
\newtheorem{remark}{Remark}[section]
\newtheorem{definition}{Definition}[section]

 {\everymath{\displaystyle\everymath{}}\array}%
 {\endarray}

\newcommand{\eps}{\varepsilon}
\renewcommand{\bar}[1]{\overline{#1}}
\newcommand{\RR}{\mathbb{R}}
\newcommand{\NN}{\mathbb{N}}
\newcommand{\Z}{\mathcal{Z}}
\newcommand{\ds}{\displaystyle}

\title{High Density Limit of the Stationary One Dimensional Schr\"odinger--Poisson
System\thanks{This work was supported by the ACI Nouvelles Interfaces des Math\'ematiques ACINIM 176-2004 entitled ``MOQUA\@,''
the French ministry of research, and the ACI Jeunes chercheurs JC1035 ``Mod\`eles dispersifs vectoriels pour
le transport \`a l'\'echelle nanom\'etrique.''}}

\author{Raymond El Hajj\thanks{Centre de Mathématiques, Insa de Rennes et IRMAR (UMR 6625),
20 avenue des Buttes de Co\"esmes, 35708 Rennes Cedex 07, France (raymond.el-hajj@insa-rennes.fr)}
\and Naoufel Ben Abdallah\thanks{Institut de Mathématiques de Toulouse (UMR 5219),
\'equipe MIP\@, Universit\'e Paul Sabatier Toulouse 3,
118 route de Narbonne, 31062 Toulouse Cedex 09, France
(naoufel@\allowbreak math.\allowbreak univ-toulouse.\allowbreak fr).}}
\date{}
\begin{document}

\maketitle

\begin{abstract}
The stationary one dimensional Schr\"odinger--Poisson system on a
bounded interval is considered in the limit of a small Debye length
(or small temperature). Electrons are supposed to be in a mixed
state with the Boltzmann statistics. Using various reformulations of
the system as convex minimization problems, we show that only the first
energy level is asymptotically occupied. The electrostatic potential
is shown to converge towards a boundary layer potential with a
profile computed by means of a half space Schr\"odinger--Poisson
system.
\end{abstract}


\textbf{AMS.} 35A15, 35J10, 35Q40, 46N50, 75G65, 81Q10
%

\section{Introduction and main results}

\subsection{Introduction}

The Schr\"odinger--Poisson system is one of the most used models for
quantum transport of charged particles in semiconductors as well as
for quantum chemistry problems
\cite{nba,BM,Castella,IZL,KaiserRehberg,weierstrass,Mark,MRS,Nier1,Nier2,Nier3,Nier4,V}.
It describes the quantum motion of an ensemble of  electrons
submitted to and interacting with an electrostatic potential. The
electron ensemble might be completely confined or in interaction
with reservoirs. In the latter case, one speaks about open systems
for which the particles are described by means of the scattering
states of the Schr\"odinger Hamiltonian corresponding to the
electrostatic potential which is in turn coupled to electron
particle density through the electrostatic interaction. This leads
to nonlinear partial differential equations whose analysis involves
scattering theory techniques and limiting absorption theorems
\cite{Nier4,nba-degmar,nba} and in which the repulsive character of
the electrostatic interaction plays an important role in the
analysis (it provides the necessary a~priori estimates for solving
the problem).

For closed systems, the particles are described thanks to the
eigenstates and eigenenergies of the Schr\"odinger Hamiltonian. The
electron density is the superposition of the densities of the
eigenstates with an occupation number decreasingly depending on
their eigenenergy. The coupling is again obtained through the
Poisson equation modeling the electrostatic interaction. This
problem was reformulated by Nier \cite{Nier1,Nier2,Nier3} as a
minimization of a convex function (whose unknown is the
electrostatic potential) which allows us to prove existence and
uniqueness results. In~\cite{weierstrass}, one can find
generalizations including local contributions to the potential and
which can be included in the functional to be minimized. This short
review partially covers stationary problems. For evolution problems,
an extended bibliography is available, and we refer the reader to the books of
Markowich, Ringhofer, and Schmeiser~\cite{MRS} and Cazenave~\cite{CAZ} for references.

In this paper we are interested in a singularly perturbed version of the Schr\"o\-dinger--Poisson 
system which arises from the description
of the so-called two dimensional electron gases
\cite{AFS,HFB}. The electrons, in such systems, are strongly
confined in one direction, at the interface between two material,
and are free to move in the two remaining ones. In~\cite{BMO}, the
analysis of the Schr\"odinger equation of strongly confined
electrons in one direction is performed. The confined direction is
called $z$ and the confining potential is assumed to be given and
scaled as $\frac{1}{\eps^2}V_c(\frac{z}{\eps})$, where $\eps$ is
a small parameter. Approximate models for the transport direction
(orthogonal to $z$) derived heuristically in the previous works
\cite{PN,nbapolizzi,P} are then analyzed in \cite{BMO,Pin}. The aim of
the present work is to somehow justify the scaling
$\frac{1}{\eps^2}V_c(\frac{z}{\eps})$ by the analysis of the
self-consistent Schr\"odinger--Poisson system in the $z$ direction. This
is why we shall forget about the transport issues in the orthogonal
direction and assume that the considered system is invariant with
respect to it. The parameter $\eps$ in the present work is linked to
the scaled Debye length as shall be explained later. The analysis
relies on the minimization formulation of the problem leading to a
singularly perturbed functional. After a rescaling argument, we are
led to the analysis of a half space Schr\"odinger--Poisson system
in which only the first eigenstate is occupied. Additional estimates
are  obtained thanks to reformulation of the single state
Schr\"odinger--Poisson system as another minimization problem whose
unknown is the first eigenfunction (and not the potential). This
formulation is used in quantum chemistry~\cite{CBL}.

Let us now come to the precise description of the problem and the results.
 The system is one dimensional and occupies the interval $[0,1]$. The
electrostatic energy is given by $V_\eps(z)$. It satisfies the
following one dimensional stationary Schr\"odinger--Poisson
system:\begin{equation}\label{mod1}
{\left\{
\begin{array}{l@{}}
\ds -\frac{d^2\varphi_{p}}{dz^2}+V\varphi_{p}= \mathcal{E}_{p}\varphi_{p}, \qquad z\in [0,1], \\ \noalign{\vspace*{\jot}}
\ds  \varphi_{p}\in H^1(0,1),\quad \varphi_{p}(0)=0,\quad \varphi_{p}(1)=0,\quad \int_0^1\varphi_{p}\varphi_{q}= \delta_{pq}, \\ \noalign{\vspace*{\jot}}
\ds -\varepsilon^3\frac{d^2V}{dz^2}= {1\over \Z} \sum_{p=1}^{+\infty}e^{-\mathcal{E}_{p}} |\varphi_{p}|^2,\quad \Z = \sum_{p=1}^{+\infty}e^{-\mathcal{E}_{p}}, \\ \noalign{\vspace*{\jot}}
\ds  V(0)=0,\quad {dV\over dz}(1)=0.
\end{array}
\right.}\hspace*{-\nulldelimiterspace}
\end{equation}
The dimensionless parameter $\varepsilon$ is a small parameter
which is devoted to tending to zero.
The choice of the third power is done for notational convenience as
shall be understood later. This parameter is related to the Debye
length and shall be explicitly given by the rescaling of the
Schr\"odinger--Poisson system (\ref{equation16}) (see subsection~1.3). The eigenvalues of the Schr\"odinger operator
$(\mathcal{E}_{p})_p$ are the energy levels in the potential well.
The sum in the right-hand side of the Poisson equation includes all
eigenvalues of the Schr\"odinger operator. In the limit $\eps\to 0$,
one expects that the wave functions concentrate at $z=0$. The
boundary condition for the potential at $z = 1$ is physically
justified in some physical situations such as in bulk materials.
However, a Dirichlet condition is more commonly used in such
problems. The analysis can be carried out in that case with the cost of
technical complexity since a new boundary layer at $z=1$
will appear and the eigenvalues will have asymptotically a double
multiplicity. For simplicity, we do not consider this case. Since
the density is very high in the limit $\eps\rightarrow 0^+$, the
Boltzmann statistics should be replaced by the Fermi--Dirac ones. The
analysis can be done in this case with the cost of technical
complications. More detailed comments about this are given in the last
section of this paper. In order to analyze the boundary layer, we make the change of
variables\begin{equation}
\label{equation1}
\varphi_{p}(z)=\frac{1}{\sqrt{\varepsilon}} \psi_{p}\left(\frac{z}{\varepsilon}\right),\quad
\mathcal{E}_{p} =\frac{1}{\varepsilon^2}E_{p}, \quad
V(z)= \frac{1}{\varepsilon^2}U\left(\frac{z}{\varepsilon}\right),\quad
\xi=\frac{z}{\eps}.
 \end{equation}
Then, $U$ verifies
$-\frac{d^2U}{d\xi^2}={1\over\tilde{\Z}}\sum_{p=1}^{+\infty} e^{\frac{-E_{p}}{\varepsilon^2}}|\psi_{p}|^2$ with
$\tilde{\Z} =\sum_{p=1}^{+\infty}e^{\frac{-E_{p}}{\varepsilon^2}}$. Since there
is a uniform gap with respect to $\eps$ between $E_1$ and $E_p$ for
$p\geq 2$ (see Lemma~\ref{gap}), the terms
$e^{\frac{-E_{p}}{\varepsilon^2}}$ with $p \geq 2$ are expected to
be negligible when compared to the first one ($p= 1$). Therefore, it
is natural to expect the solution of \eqref{mod1} to be
asymptotically close to the solution of the following
Schr\"odinger--Poisson system in which only the first energy level is taken into
account:\begin{equation}
\label{mod2}
{\left\{
\begin{array}{l@{}}
\ds -\frac{d^2\tilde\varphi_{1}}{dz^2}+ \tilde V\tilde\varphi_{1}= \tilde{\mathcal{E}}_{1}\tilde\varphi_{1}, \qquad z\in [0,1],\\ \noalign{\vspace*{\jot}}
\ds \tilde{\mathcal{E}}_{1}=\inf_{\varphi\in H^1_0(0,1),\, \|\varphi\|_{L^2}=1}\left\{\int_0^1|\varphi^\prime|^2 +\int_0^1\tilde V\varphi^2\right\}, \\ \noalign{\vspace*{\jot}}
\ds -\varepsilon^3\frac{d^2\tilde V}{dz^2}= |\tilde \varphi_{1}|^2,\\ \noalign{\vspace*{\jot}}
\ds \tilde V(0)=0,\quad {d \tilde V \over dz}(1)=0.
\end{array}
\right.}\hspace*{-\nulldelimiterspace}
\end{equation}
Moreover, when $\eps$ goes to zero, we will prove that the
electrostatic potential, $\tilde V_\varepsilon$, solution of
\eqref{mod2} converges towards a boundary layer potential with
profile, $U_0$, solution of the following half line
problem:\begin{equation}
\label{mod5}
{\left\{
\begin{array}{l@{}}
\ds -\frac{d^2\psi_{1}}{d\xi^2}+U\psi_{1}=E_{1}\psi_{1},\qquad \xi\in [0,+\infty[, \\  \noalign{\vspace*{\jot}}
\ds E_{1}=\inf_{\psi\in H^1_0(\mathbb{R}^+),\, \|\psi\|_{L^2}=1}\left\{\int_0^{+\infty}|\psi^\prime|^2 +\int_0^{+\infty}U\psi^2\right\},\\ \noalign{\vspace*{\jot}}
\ds -\frac{d^2U}{d\xi^2}=|\psi_{1}|^2, \\ \noalign{\vspace*{\jot}}
\ds U(0)=0,\quad {d U\over d\xi}\in L^2(\mathbb{R}^+).
\end{array}
\right.}\hspace*{-\nulldelimiterspace}
\end{equation}

\subsection{Main results}

In this paper, a rigourous analysis and comparison of the systems presented above
will be provided. Namely, \eqref{mod1} and \eqref{mod2} are posed on a bounded domain.
The one dimensional Schr\"odinger--Poisson system on a bounded
interval was studied by Nier in~\cite{Nier1}. Each of these systems can be reformulated
as a minimization problem (see section~2 for details). However, the limit problem \eqref{mod5} is posed
on an unbounded domain. Our first result deals with the study of \eqref{mod5}. We also prove that
it can be formulated as a minimization problem.

\begin{theorem}\label{mainth1}
Let $J_0(.)$ be the energy functional
defined on $\dot H^1_0(\mathbb{R}^+)$ (given by\/ {\rm(\ref{equation26})})
by\begin{equation}
\label{equation31}
J_0(U)=\frac{1}{2}\int_0^{+\infty}|U^\prime|^2-E_1^\infty[U],
\end{equation}
where $E_1^\infty[U]$ is the fundamental mode of the
Schr\"odinger operator given by\/ {\rm(\ref{equation36})}. The limit
problem\/ {\rm(\ref{mod5})} has a unique solution\/ $(U_0, E_{1,0},\psi_{1,0})$, and $U_0$ satisfies the following minimization
problem:\begin{equation}
\label{equation35}
J_0(U_0)=\inf_{U\in \dot{H}^1_0(\mathbb{R}^+)}J_0(U).
\end{equation}\end{theorem}\unskip

The comparison of the systems presented above is established by our second main theorem.

\begin{theorem}\label{mainth2}
Let $V_\eps$, $\tilde V_\eps$, and $U_0$ be the potentials
satisfying problems\/ {\rm(\ref{mod1})},\/ {\rm(\ref{mod2})}, and\/ {\rm(\ref{mod5})}, respectively. Then the following estimates
hold:\begin{equation}\label{equation18}
\| V_\varepsilon-\tilde V_\eps\|_{H^1(0,1)}=\mathcal{O}(e^{-\frac{c}{\varepsilon^2}})
\end{equation}
and\begin{equation}
\label{equation17}
\left\|\tilde V_\varepsilon- \frac{1}{\eps^2}U_0\left(\frac{.}{\eps}\right)\right\|_{H^1(0,1)}= \mathcal{O}(e^{-\frac{c}{\varepsilon}}),
\end{equation}
where $c$ is a general strictly positive constant
independent of~$\varepsilon$.\end{theorem}

The paper is organized as follows. In the next subsection, we
present some remarks on the scaling giving model \eqref{mod1}, and we
end this section by fixing some notation and definitions. In
section~\ref{S-Pbounded}, we recall the spectral properties of the
Schr\"odinger operator on a bounded domain and state the optimization
problems corresponding to \eqref{mod1} and \eqref{mod2} (or more
precisely to the intermediate systems \eqref{mod3} and
\eqref{mod4}). Section~\ref{limit} is devoted to the analysis of the
limit problem (\ref{mod5}) posed on the half line (proof of
Theorem~\ref{mainth1}). We will first study the properties of the
fundamental mode of the Schr\"odinger operator (Proposition~\ref{proposition3.1}). The limit problem leads us to the study of a
minimization problem posed on an unbounded domain. This will be done by
means of the concentration-compactness principle introduced by Lions
in~\cite{lions}. Estimates (\ref{equation18}) and (\ref{equation17})
are proved in section~\ref{convergence}. Some comments concerning
the Fermi--Dirac statistics, the choice of the boundary conditions,
and the problems of the multidimensional case are given in
section~\ref{comments}. Finally, Appendix~\ref{appendix} is devoted to the
proof of Lemma~\ref{lemmaA.1}.

First, let us make this remark.

\begin{remark}
To prove\/ {\rm(\ref{equation18})--(\ref{equation17})}, we use the scaled versions of\/ {\rm(\ref{mod1})} and\/ {\rm(\ref{mod2})}
when applying the changes of variables\/ {\rm(\ref{equation1})}
and\begin{equation}
\label{equation1'}
\tilde \varphi_{1}(z)=\frac{1}{\sqrt{\varepsilon}} \tilde \psi_{1}\left(\frac{z}{\varepsilon}\right),\quad
\tilde {\mathcal{E}}_{1} =\frac{1}{\varepsilon^2}\tilde E_{1}, \quad
\tilde V(z)=  \frac{1}{\varepsilon^2}\tilde U\left(\frac{z}{\varepsilon}\right),\quad
\xi=\frac{z}{\eps}.
\end{equation}
Then the intermediate Schr\"odinger--Poisson models
write\begin{equation}
\label{mod3}
{\left\{
\begin{array}{l@{}}
\ds -\frac{d^2\psi_{p}}{d\xi^2}+U\psi_{p}= E_{p}\psi_{p}, \qquad \xi \in \left[0,\frac{1}{\varepsilon}\right],\\ \noalign{\vspace*{\jot}}
\ds \psi_{p}\in H^1\left(0,\frac{1}{\varepsilon}\right),\quad \psi_{p}(0)=0,\quad \psi_{p}\left(\frac{1}{\eps}\right)=0,\quad
\int_0^{\frac{1}{\varepsilon}} \psi_{p}\psi_{q}=\delta_{pq}, \\ \noalign{\vspace*{\jot}}
\ds -\frac{d^2U}{d\xi^2}={1\over\tilde{\Z}}\sum_{p=1}^{+\infty} e^{\frac{-E_{p}}{\varepsilon^2}}|\psi_{p}|^2, \quad
\tilde{\Z} =\sum_{p=1}^{+\infty}e^{\frac{-E_{p}}{\varepsilon^2}}, \\ \noalign{\vspace*{\jot}}
\ds U(0)=0,\quad {dU\over d\xi}\left(\frac{1}{\varepsilon}\right)=0
\end{array}
\right.}\hspace*{-\nulldelimiterspace}
\end{equation}
and\begin{equation}
\label{mod4}
{\left\{
\begin{array}{l@{}}
\ds -\frac{d^2\tilde\psi_{1}}{d\xi^2}+\tilde U\tilde\psi_{1}= \tilde E_{1}\tilde\psi_{1},\qquad \xi\in \left[0,\frac{1}{\varepsilon}\right], \\ \noalign{\vspace*{\jot}}
\ds \tilde E_{1}=\inf_{\psi\in H^1_0(0,\frac{1}{\varepsilon}),\, \|\psi\|_{L^2}=1} \left\{\int_0^{\frac{1}{\varepsilon}}|\psi^\prime|^2
+\int_0^{\frac{1}{\varepsilon}}\tilde U\psi^2\right\},\\  \noalign{\vspace*{\jot}}
\ds -\frac{d^2\tilde U}{d\xi^2}=|\tilde \psi_{1}|^2, \\  \noalign{\vspace*{\jot}}
\ds \tilde U(0)=0,\quad {d\tilde U\over d\xi}\left(\frac{1}{\varepsilon}\right)=0.
\end{array}
\right.}\hspace*{-\nulldelimiterspace}
\end{equation}
Remark that it is natural to expect\/ {\rm(\ref{mod4})} to be close, when
$\eps$ goes to zero, to the limit problem\/ {\rm(\ref{mod5})} posed on\/ $[0,+\infty)$.\end{remark}

\subsection{Remark on the scaling}

Here we show how the system (\ref{mod1}) can be obtained by a
rescaling of the Schr\"odinger--Poisson system written with the
physical dimensional variables. Indeed, let $(\chi_p(Z),\Lambda_p)$
be the eigenfunctions and the eigenenergies of the one dimensional
Schr\"odinger operator (the confinement operator)
$-\frac{\hbar^2}{2m}\frac{d^2}{dZ^2}+W$ with
homogeneous Dirichlet
data:\begin{equation}
\label{equation8}
-\frac{\hbar^2}{2m}\frac{d^2\chi_p}{dZ^2}+W\chi_p=\Lambda_p\chi_p,
\end{equation}
where $\hbar$ is the Planck constant and $m$ denotes the
effective mass of the electrons in the crystal. The $(\chi_p)_p$ is
an orthonormal basis of $L^2(0,L)$. The variable $Z$ belongs to
$[0,L]$, where $L$ is the typical length of the confinement. Denoting
by $n$ the electronic density, this can be
written\begin{equation}
\label{equation9}
n(Z)=\sum_{p=1}^{+\infty}n_p|\chi_p(Z)|^2.
\end{equation}
In this formula, $|\chi_p(Z)|^2$ is the probability of presence at point $Z$
of an electron in the $p$th state. Using Boltzmann statistics, the occupation factor $n_p$ is given
by\begin{equation}
\label{equation10}
n_p= {N_s \over {\cal Z}} \exp\left(-\frac{\Lambda_p}{k_BT}\right), \quad
{\cal Z}= \sum_{q=1}^{+\infty} \exp\left(-\frac{\Lambda_q}{k_BT}\right),
\end{equation}
where $k_B$ is the Boltzmann constant, $T$ denotes the
temperature, and $N_s$ is the surface density assumed to be given.
With this notation we have $\int_0^L n(Z)\,  dZ = N_s$, which
means that the total number of electrons in the interval $[0,L]$
(per unit surface in the two remaining  spatial directions) is
given.
 The electrostatic potential $W$ and the electron density $n$ are coupled through the Poisson
 equation:\begin{equation}
 \label{equation11}
-\frac{d^2W}{dZ^2}=\frac{q^2}{\varepsilon_0\varepsilon_r}n
\end{equation}
with boundary
conditions\begin{equation}
\label{equation12}
W(0)=0,\quad
\frac{dW}{dZ}(L)=0.
\end{equation}
In (\ref{equation11}), the constant $q$ is the elementary
electric charge and $\varepsilon_0$, $\varepsilon_r$ are,
respectively, the permittivity of the vacuum and the relative
permittivity of the material.

Let us rescale the problem (\ref{equation8})--(\ref{equation12}) by noticing
that\begin{equation}
\label{equation13}
z=\frac{Z}{L}\in[0,1],\quad
W(Z)= (k_BT)V\left(\frac{Z}{L}\right), \quad
\Lambda_p=(k_BT)\mathcal{E}_{p}, \quad
\chi_p(Z)=\frac{1}{\sqrt{L}}\varphi_{p}\left(\frac{Z}{L}\right).
\end{equation}
We assume that $\frac{\hbar^2}{2mL^2}$
is of the same order of the thermal energy $(k_BT)$.
In order to simplify the mathematical presentation, we suppose
that\begin{equation}
\label{equation14}
\frac{\hbar^2}{2mL^2}=k_BT.
\end{equation}
By inserting (\ref{equation13}) into the system (\ref{equation8})--(\ref{equation12}),
we obtain, after straightforward computation, the system (\ref{mod1}) in which $\varepsilon$ is related to the scaled Debye
length:\begin{equation}
\label{equation16}
\varepsilon^3 = \left(\frac{\lambda_D}{L}\right)^2,\quad
\lambda_D = \sqrt{\frac{k_BT\varepsilon_0\varepsilon_r}{q^2N}},
\end{equation}
where $N = {N_s\over L}$ is the average volume density of electrons.

\subsection{Notation and definitions}

We summarize in this subsection the different variables and
notation used in this paper.
\begin{itemize}
\item For the
Schr\"odinger--Poisson problems posed on $[0,1]$, $z$ denotes the
space variable, $V$ denotes the potential variable, and
$(\mathcal{E},\varphi)$ represents any eigenvalue and the
corresponding eigenfunction of the Schr\"odinger operator. For
systems posed on $[0,\frac{1}{\eps}]$ or on $\RR^+$, we use
$\xi$, $U$, and $(E,\psi)$ as variables. The same notation with
$\tilde{~}$, i.e., $(\tilde V,\tilde{\mathcal{E}},\tilde\varphi)$ or
$(\tilde U,\tilde E,\tilde\psi)$, is used for the variables of
Schr\"odinger--Poisson systems in which only the first eigenstate is
taken into account.
\item For any real valued function $V\in L^2(0,L)$, where $L >0$
is given ($L=1$ or $\frac{1}{\eps}$ here), we denote by
$H[V]$ the Dirichlet--Schr\"odinger
operator\begin{equation}
\label{equation19}
H[V]=-\frac{d^2}{dx^2}+V(x) \quad (x =z \text{ or } \xi \text{ here})
\end{equation}
defined on the domain $D(H[V])=H^2(0,L)\cap H^1_0(0,L)$.
In addition, the sequence of eigenenergies and eigenfunctions of
$H[V]$ will be denoted by $(E_p[V],\psi_p[V])_{p\in \NN^*}$. We give
in the next section the main properties satisfied by the functions
$V\mapsto E_p[V]$ and $V\mapsto \psi_p[V]$ for any $p\in \NN^*$.
\item The potentials satisfying \eqref{mod1} and \eqref{mod2}
are denoted by $V_\eps$ and $\tilde V_\eps$. In addition,
$(\mathcal{E}_{p,\eps},\varphi_{p,\eps})$ and
$(\tilde{\mathcal{E}}_{p,\eps},\tilde{\varphi}_{p,\eps})$, with
$p\in \NN^*$, represent the corresponding energy couples of
$H[V_\eps]$ and $H[\tilde V_\eps]$, respectively. In other words,
$\mathcal{E}_{p,\eps}:=E_p[V_\eps]$,
$\varphi_{p,\eps}:=\psi_p[V_\eps]$,
$\tilde{\mathcal{E}}_{p,\eps}:=E_p[\tilde V_\eps]$, and
$\tilde{\varphi}_{p,\eps}:=\psi_p[\tilde V_\eps]$. Similarly, the
solutions of \eqref{mod3} and \eqref{mod4} will be denoted,
respectively, by $(U_\eps,E_{p,\eps},\psi_{p,\eps})$ and $(\tilde
U_\eps,\tilde E_{p,\eps},\tilde{\psi}_{p,\eps})$. Finally, we fix
$(U_0,E_{1,0},\psi_{1,0})$ to denote the solution of the limit
problem \eqref{mod5}.
\end{itemize}

Let us now define some spaces which will be
used throughout this paper.

\begin{definition}\label{definition2.1}
{\rm(i)} For $L>0$, we  define
\begin{equation}\label{equation25}
H^{1,0}(0,L)= \left\{U\in {H^1(0,L),}~U(0)=0\right\}.
\end{equation}

{\rm(ii)}  The space $\dot{H}^1_0(\mathbb{R}^+)$ is defined as
follows:\begin{equation}
\label{equation26}
\dot{H}^1_0(\mathbb{R}^+)= \bigl\{U\in {L^2_{loc}(\mathbb{R}^+),}\text{ } U^\prime \in {L^2(\mathbb{R}^+),}\text{ }U(0)={0,}\text{ and } U\geq 0\bigr\}\,.
\end{equation}

{\rm(iii)} For any\/ $0<L\leq +\infty$, we shall denote by $S_L$ the set
of normalized functions of $H^1_0(0,L)$ with respect to the
$L^2$-norm\begin{equation}
\label{equation27}
S_L= \left\{\varphi \in {H^1_0(0,L),}~ {\int_0^{L}\varphi^2} =1\right\}.
\end{equation}
Here $H^1_0(0,L)$ is the space of $H^1$-functions
vanishing on\/ {\rm0}~and~$L$, and when $L=+\infty$
\[
H^1_0(\RR^+)=\left\{\psi\in {H^1(\RR^+),}\text{ }\psi(0)= {0} \right\}.
\]\end{definition}\unskip

\section{Schr\"odinger--Poisson system on a bounded domain}\label{S-Pbounded}

We begin this part by recalling some basic
properties satisfied by the eigenvalues and the eigenfunctions of
the one dimensional Schr\"odinger operator \eqref{equation19}. These
properties are standard and can be found in
\cite{Kato,Nier1,Poschel,Reed Simon}. The operator $H[V]$ is
self-adjoint, is bounded from below, and has compact resolvent. There
exists a strictly increasing sequence $(E_p[V])_p$ of real numbers
tending to $+\infty$ and an orthonormal basis of $L^2(0,L)$,
$(\psi_p[V])_p$, such that $\psi_p[V]\in D(H[V])$
and\begin{equation}
\label{equation21}
H[V]\psi_p[V]=E_p[V]\psi_p[V].
\end{equation}
For $V=0$, we have by a simple
calculation\begin{equation}
\label{equation21'}
E_p[0]=\frac{\pi^2p^2}{L^2},\quad
\psi_p[0](x)= \sqrt{\frac{2}{L}}\sin\left(\frac{p \pi x}{L}\right).
\end{equation}
The eigenvalues $E_p[V]$ are simple and satisfy the following characterization (min-max
principle)~\cite{Reed Simon}:\begin{equation}
\label{equation22}
E_p[V]=\min_{V_p\in \mathbb{V}_p\left(D(H[V])\right)}\max_{\varphi\in V_p,\, \varphi\neq 0}\frac{(H[V]\varphi,\varphi)_{L^2}}{\|\varphi\|_{L^2}^2},
\end{equation}
where $\mathbb{V}_p\left(D(H[V])\right)$ is the set of the
subspaces of $D(H[V])$ with dimension equal to $p$, and $(.,.)$ denotes
the scalar product in $L^2$. In view of the min-max formula
(\ref{equation22}), one can verify that for any $p\in\mathbb{N}^*$,
$E_p[.]$ is an increasing function, which means
that\[
E_p[V]\leq E_p[W]\quad
\text{if } V\leq W \text{ a.e.}
\]
Moreover, we have the Lipschitz
property, for any real valued functions $V,W$ in
$L^\infty(0,L)$,\begin{equation}
\label{equation23}
|E_p[V]-E_p[W]|\leq \|V-W\|_{L^\infty(0,L)}.
\end{equation}
Besides, one can prove the following lemma~\cite{Nier1}.

\begin{lemma}
For any $p\in \mathbb{N}^*$, the
maps\[
E_p[.]:L^2(0,L)\longrightarrow \mathbb{R}, \quad
\psi_p[.]:L^\infty(0,L)\longrightarrow L^1(0,L)
\]
are G\^ateaux differentiable, and their derivatives are given, respectively,
by\begin{equation}
\label{equation23'}
\begin{split}
dE_p[V].W &= \int_0^{L}|\psi_p[V]|^2W dx\quad\text{and} \\
d\psi_p[V].W &=\sum_{q\neq p}\frac{1}{E_p[V]-E_q[V]}\left(\int_0^{L}W\psi_p\psi_q dx\right)\psi_q
\end{split}
\end{equation}
for any $V,W \in L^\infty(0,L)$.\end{lemma}

Using the spectral properties of the Schr\"odinger
operator, one can prove the following proposition. For details on
the proof see~\cite{Nier1}.

\begin{proposition}
The systems\/ {\rm(\ref{mod3})} and\/ {\rm(\ref{mod4})} are well posed.
They are equivalent, respectively, to the following minimization
problems:\begin{equation}
\label{equation33}
J_\varepsilon(U_\varepsilon)=\inf_{U\in H^{1,0}(0,M_\varepsilon)}J_\varepsilon(U)
\end{equation}
and\begin{equation}
\label{equation34}
\tilde J_\varepsilon(\tilde U_\varepsilon)=\inf_{U\in H^{1,0}(0,M_\varepsilon) } \tilde J_\varepsilon(U),
\end{equation}
where $M_\varepsilon=\frac{1}{\varepsilon}$. The
energy functionals $J_\varepsilon$ and $\tilde J_\varepsilon$ are
given by\begin{equation}
\label{equation29}
J_\varepsilon(U)=\frac{1}{2}\int_0^{M_\varepsilon}|U^\prime|^2+ \varepsilon^2 \log\left(\,\sum_{p=1}^{+\infty}e^{-\frac{E_p[U]}{\varepsilon^2}}\right)
\end{equation}
and\begin{equation}
\label{equation29'}
\tilde J_\varepsilon(U)=\frac{1}{2}\int_0^{M_\varepsilon}|U^\prime|^2-E_1[U].
\end{equation}
Each one of problems\/ {\rm(\ref{equation33})} and\/ {\rm(\ref{equation34})} admits a unique solution.\end{proposition}

\begin{remark}
One can similarly study the systems\/ {\rm(\ref{mod1})} and\/ {\rm(\ref{mod2})}
and prove that each one is equivalent to an optimization problem.\end{remark}

\section{Analysis of the limit problem (\ref{mod5})}\label{limit}

The aim of this part is to study the well-posedness of the limit
problem (\ref{mod5}) posed on the half line. Namely, this part is
concerned with the proof of Theorem~\ref{mainth1}. We begin with the
study of the fundamental mode, $E_1^\infty[.]$, of the Schr\"odinger
operator. Its main properties are listed in Proposition~\ref{proposition3.1}.

\subsection{Properties of the fundamental mode of the Schr\"odinger operator on \boldmath$[0,+\infty)$}

We begin by defining the fundamental mode.

\begin{definition}
For any real and positive function $U\in L^1_{loc}(\RR^+) $, the
fundamental mode of the Schr\"odinger operator
is\begin{equation}
\label{equation36}
E_1^\infty[U]=\inf_{\psi\in S_\infty}J_U(\psi),
\end{equation}
where for any $\psi\in S_\infty$ (defined by\/ {\rm(\ref{equation27})}) we
have\begin{equation}
\label{equation36'}
J_U(\psi)=\int_0^{+\infty}|\psi^\prime|^2+\int_0^{+\infty}U\psi^2.
\end{equation}\end{definition}\unskip

One difficulty due to the unboundedness of the interval $[0,+\infty)$ is that $E_1^\infty[.]$
might not be an eigenvalue but only the lower bound of
the essential spectrum. The following proposition gives
some properties of $E_1^\infty[.]$ and some sufficient conditions on the potential
for which $E_1^\infty[.]$ is an eigenvalue.

\begin{proposition}\label{proposition3.1}
\begin{enumerate}
\item The map $U\mapsto E_1^\infty[U]$
is a continuous, concave, and  increasing function with values in\/ $\overline{\RR^+}:=[0,+\infty]$
satisfying\begin{equation}
\label{equation37}
E_1^\infty[U]\leq \limsup_{\xi\rightarrow+\infty}U(\xi).
\end{equation}
\item If $U\in L^1_{loc}(\RR^+)$, $U\geq 0$ such that
$E_1^\infty[U]<\liminf_{\xi\rightarrow+\infty}U(\xi)$,
then $E_1^\infty[U]$ is reached by a unique positive function
$\psi_1[U]$, which means that there exists a unique positive function
$\psi_1[U]\in S_\infty$ such that $E_1^\infty[U]=J_U(\psi_1[U])$. In addition, we
have\begin{equation}
\label{equation38}
\frac{dE_1^\infty}{dU}[U].W=\int_0^{+\infty}|\psi_1[U]|^2Wd\xi
\end{equation}
for any function $W$ in $L_0^\infty(\RR^+)$, the space of
bounded functions with compact support on\/~$\RR^+$.
\item Let $U\in L^1_{loc}(\RR^+)$ be a positive function such that\/ $\lim_{\xi\rightarrow+\infty}U(\xi)$ exists,
$U\leq \lim_{\xi\rightarrow+\infty}U(\xi)$, and
$E_1^\infty[U]=\lim_{\xi\rightarrow+\infty}U(\xi)$.  Then we
have\[
\frac{dE_1^\infty}{dU}[U].W=0
\]
for any $W\in L_0^\infty(\RR^+)$.
\item Let $\alpha$ be an arbitrary positive constant.  Then we
have\begin{equation}
\label{equation38'}
E_1^\infty[\alpha\sqrt{\xi}]=\alpha^{\frac{4}{5}}E_1^\infty[\sqrt{\xi}].
\end{equation}
\end{enumerate}\end{proposition}

\begin{remark}
There is quite a difference between the third case of this
proposition, where $E_1^\infty[U]=\lim_{\xi\rightarrow
+\infty}U(\xi)$, and the second case, which includes $E_1^\infty[U]< \lim_{\xi\rightarrow +\infty}U(\xi)$. This result is
natural and can be interpreted as follows. The classically allowed
region for a particle with energy $E$ is the set $\mathcal{A}=\{\xi \in [0,+\infty);${\rm~}$U(\xi)\leq E\}$. In the case $E\geq U(\xi)$
on\/ $[0,+\infty)$, the set $\mathcal{A}$ extends to\/ $+\infty$ so that
there is no bound state, while in the case $E<\lim_{\xi\rightarrow +\infty}U(\xi) $ the set $\mathcal{A}$ is
bounded and $E$ is a bounded state energy.\end{remark}

\begin{lemma}\label{lemmaA.1}
Let $U\in \dot H^1_0(\mathbb{R^+})$ such that $E_1^\infty[U]<\liminf_{+\infty}U$.
Then all minimizing
sequences\/ $(\psi_n)_n$ of problem\/ {\rm(\ref{equation36})} are relatively compact in
$L^2(\mathbb{R}^+)$.\end{lemma}

 This lemma is needed for the proof of the second point of Proposition~\ref{proposition3.1}. It is proved in
 Appendix~\ref{appendix}\@. The
proof is based on the concentration-compactness principle.

{\it Proof of Proposition\/~{\rm\ref{proposition3.1}}}. 1.
Remark first that for any positive function $U$, $E_1^\infty[U]$
exists and belongs to $\overline{\RR^+}$. It is easy to check, from
the definition of $E_1^\infty[.]$, that it is a continuous, concave,
and increasing function. To prove inequality (\ref{equation37}), let
$\psi\in S_\infty$ be fixed and set $\psi_\delta=\sqrt{\delta}\psi(\delta \xi)$ for any real positive
$\delta$. Then $\psi_\delta\in S_\infty$, and since $E_1^\infty[U]$
verifies (\ref{equation36}), we
have\begin{equation}
\label{equation39}
E_1^\infty[U]\leq J_U(\psi_\delta).
\end{equation}
Moreover, we have $J_U(\psi_\delta)=\delta^2\int_0^{+\infty}|\psi^\prime(\xi)|^2d\xi + \int_0^{+\infty}U(\frac{\xi}{\delta})\psi^2(\xi)d\xi$.
Then\[
\limsup_{\delta\rightarrow 0}J_U(\psi_\delta)\leq \int_0^{+\infty}\limsup_{\delta\rightarrow 0}U\left(\frac{\xi}{\delta}\right)\psi^2(\xi)d\xi\leq \limsup_{\xi\rightarrow+\infty}U(\xi).
\]
Taking the $\limsup_{\delta\rightarrow 0}$
of (\ref{equation39}), one obtains inequality (\ref{equation37}).

 2. Let $(\psi_n)_n$ be a minimizing sequence of $E_1^\infty[U]$; i.e., $\psi_n\in S_\infty$ for any $n\in \mathbb{N}^*$ and
$J_U(\psi_n)\rightarrow_{n\rightarrow +\infty} E_1^\infty[U]$. The
sequence $(\psi_n)_n$ is bounded in $H^1_0(\mathbb{R}^+)$, there
exist a function $\psi\in H^1_0(\mathbb{R}^+) $ and a subsequence
also denoted $(\psi_n)$ such that $(\psi_n)$ converges weakly to
$\psi$ in $H^1_0(\mathbb{R}^+)$, and since $J_U(.)$ is weakly lower
semicontinuous (it is strictly convex and lower semicontinuous) we
have $J_U(\psi)\leq \liminf_{n\rightarrow +\infty} J_U(\psi_n)$.
Then\begin{equation}
\label{equation40}
J_U(\psi)\leq E_1^\infty[U].
\end{equation}
Besides, the hypothesis $E_1^\infty[U]<\liminf_{+\infty}U$ implies that the sequence
$(\psi_n)_n$ is relatively compact in $L^2(\mathbb{R}^+)$ (see Lemma~\ref{lemmaA.1}). Then, up to an extraction of subsequence,
$(\psi_n)_n$ converges strongly to $\psi$ in $L^2(\mathbb{R}^+)$.
Since $\|\psi_n\|^2_{L^2(\RR^+)}=1$, for all $n$, we have
$\|\psi\|^2_{L^2(\RR^+)}=1$, and then $\psi$ belongs to $S_\infty$.
Therefore, in view of the definition of $E_1^\infty[U]$
(\ref{equation36}), $E_1^\infty[U]\leq J_U(\psi)$ and with
(\ref{equation40}) we have $E_1^\infty[U]=J_U(\psi)$. Let us now
show that $E_1^\infty[U]$ is a simple eigenvalue and the
corresponding eigenfunction has a constant sign. Indeed, let
$\psi_1$ and $\psi_2$  be two minimizers of $J_U(.)$ on $S_\infty$,
i.e., $\psi_1, \psi_2\in S_\infty$, such that $E_1^\infty[U]=J_U(\psi_1)=J_U(\psi_2)$,  and let
$\phi=\sqrt{\frac{\psi_1^2}{2}+\frac{\psi_2^2}{2}}$. The function $\phi$ belongs to $S_\infty$, and we
have\[
J_U(\phi)=\frac{1}{2}J_U(\psi_1)+\frac{1}{2}J_U(\psi_2)- \int_0^{+\infty}\left|\frac{\psi_1\psi_2'-\psi_2\psi_1'}{2\phi}\right|^2
=E_1^\infty[U]-\int_0^{+\infty}\left|\frac{\psi_1\psi_2'-\psi_2\psi_1'}{2\phi}\right|^2.
\]
Since $E_1^\infty[U]\leq J_U(\phi)$  ($\phi\in S_\infty$),
we get $\int_0^{+\infty}\big|\frac{\psi_1\psi_2'-\psi_2\psi_1'}{2\phi}\big|^2=0$,
which implies that $\psi_1$ and $\psi_2$ are proportional, and so
$E_1^\infty[U]$ is simple. In particular, $\psi$ and $|\psi|$ are
two minimizers of $E_1^\infty[U]$; they are then proportional, and
since $\int_0^{+\infty}|\psi|^2=1$ we conclude that
$\psi=\pm|\psi|$. We then choose $\psi_1[U]=|\psi|$, which is
positive. This is the unique positive eigenfunction corresponding to
$E_1^\infty[U]$. To end the proof of the second point of Proposition~\ref{proposition3.1}, let $W$ be a compactly supported bounded
function ($W\in L_0^\infty(\mathbb{R}^+)$) and remark that for a
small real $t$ we have
$E_1^\infty[U+tW]\leq E_1^\infty[U]+ |t|\|W\|_\infty< \liminf_{\xi\rightarrow+\infty}U(\xi)$. In addition, since $W\in L_0^\infty$, we have
$\liminf_{\xi\rightarrow+\infty}U=\liminf_{\xi\rightarrow+\infty}(U(\xi)+tW(\xi))$.
Then, for any small real $t$, we have $E_1^\infty[U+tW]<\lim_{\xi\rightarrow+\infty}(U(\xi)+tW(\xi))$.
Therefore, for all bounded and compactly supported functions $W$ and
for all $t\in \mathbb{R}$ small, $E_1^\infty[U+tW]$ is an
eigenvalue. Let $\psi_t$ be the corresponding positive eigenfunction. We
have\[
E_1^\infty[U+tW]=\int_0^{+\infty}|\psi_t^\prime|^2d\xi+\int_0^{+\infty}\,(U+tW)|\psi_t|^2d\xi \geq E_1^\infty[U]+t\int_0^{+\infty} W|\psi_t|^2d\xi.
\]
Similarly, one
has\[
E_1^\infty[U]\geq
E_1^\infty[U+tW]-t\int_0^{+\infty}W|\psi_1[U]|^2d\xi.
\]
Then, if $t$ is a small nonnegative real (without loss of generality), one can
write\begin{equation}
\label{equation41}
\int_0^{+\infty}W|\psi_t|^2d\xi \leq
\frac{E_1^\infty[U+tW]-E_1^\infty[U]}{t}\leq
\int_0^{+\infty}W|\psi_1[U]|^2d\xi.
\end{equation}
Besides, since $(\psi_t)_t$ is bounded in
$H^1_0(\mathbb{R}^+)$, there exists a positive function $\psi_0\in
H^1_0$ such that $\psi_t$ converges weakly to $\psi_0$, when
$t\rightarrow 0^+$, in $H^1_{loc}(\mathbb{R}^+)$ and strongly in
$L^2_{loc}(\mathbb{R}^+)$. By passing to the limit $t\rightarrow 0^+$
in\[
-\psi_t''+(U+tW)\psi_t= E_1^\infty[U+tW]\psi_t
\]
we
obtain\[
-\psi_0''+U\psi_0=E_1^\infty[U]\psi_0\quad \text{in }
\mathcal{D}^\prime(0,+\infty).
\]
Since $\psi_0$ is positive, we deduce that $\psi_0=\psi_1[U]$.
Finally, to obtain (\ref{equation38}) we just have to take the limit
$t\rightarrow 0^+$ of (\ref{equation41}).

3. Remark first that, since $E_1^\infty[.]$ is a  nondecreasing real function, we have for
$t\geq 0$\[
\frac{E_1^\infty[U-t|W|]-E_1^\infty[U]}{t}\leq
\frac{E_1^\infty[U+tW]-E_1^\infty[U]}{t}\leq
\frac{E_1^\infty[U+t|W|]-E_1^\infty[U]}{t}.
\]
Therefore, it is sufficient to prove ${dE_1^\infty\over dU}[U].W
=0$ for  $W\geq 0$
 and $W\leq 0$ (the general case can be deduced by
 passing to the limit $t\rightarrow 0^+$ in the above inequalities).

(i) Let $W\in L_0^\infty(\RR^+)$ and $W\geq 0$.  Then we have
$E_1^\infty[U]\leq E_1^\infty[U+tW]$. Besides, by
(\ref{equation37}), we have $E_1^\infty[U+tW]\leq
\lim_{\xi\rightarrow+\infty}(U(\xi)+tW(\xi))=\lim_{\xi\rightarrow+\infty}U(\xi)=E_1^\infty[U]$.
Then, for all $W\geq 0$ in $L_0^\infty$,
$E_1^\infty[U+tW]=E_1^\infty[U]$, and the result is proved in this
case.

(ii) Let $W\in L_0^\infty(\RR^+)$, let
$W\leq 0$, and let $(t_n)_{n\in \mathbb{N}}$ be a sequence decreasing
towards~$0^+$. The sequence $(E_1^\infty[U+t_nW])_n$ is increasing
and satisfies $E_1^\infty[U+t_nW]\leq \lim_{+\infty}U
= E_1^\infty[U]$ for all $n\in \mathbb{N}$. Therefore, either it is
stationary in the vicinity of $+\infty$ and in that case
${dE_1^\infty\over dU}[U].W =0$, or it
satisfies\begin{equation}
\label{equation42}
E_1^\infty[U+t_nW]<\lim_{+\infty}U \quad \forall  n \in \mathbb{N}.
\end{equation}
In the latter case, $E_1^\infty[U+t_n W]$ is an eigenvalue and there
exists a sequence
 $(\psi_n)\in S_\infty$, $\psi_n\geq 0$, such
 that\[
E_1^\infty[U+t_nW]=J_{U+t_nW}(\psi_n)= \inf_{\psi\in
S_\infty}J_{U+t_n W}(\psi).
\]
Besides, we have $E_1^\infty[U+t_nW]\geq E_1^\infty[U] + t_n \int_0^{+\infty}W\psi_n^2d\xi$
and\begin{equation}
\label{equation43}
\left|\frac{E_1^\infty[U+t_nW]-E_1^\infty[U]}{t_n} \right|\leq
-\int_0^{+\infty}W\psi_n^2d\xi.
\end{equation}
The sequence $(\psi_n)_n$ being bounded in
$H^1_0(\mathbb{R^+})$, one can find a positive function $\psi\in
H^1_0(\mathbb{R}^+)$ and a subsequence of $(\psi_n)_n$ also denoted
by $(\psi_n)_n$ such that $\psi_n$ converges weakly to $\psi$ in
$H^1_{loc}(\mathbb{R}^+)$ and strongly in $L^2_{loc}(\mathbb{R}^+)$.
In addition $\psi$ satisfies, in the sense of
distributions,\[
-\psi''+U\psi=E_1^\infty[U]\psi =\lim_{+\infty}(U) \psi.
\]
This implies that $\psi''= (U-\lim_{+\infty}U)\psi\leq
0$ with $\psi \in H^1_0(\mathbb{R}^+)$. We  deduce that $\psi=0$ a.e., and we get the result by passing to the limit in
(\ref{equation43}), $W$ being compactly supported.

4. Let us now verify the identity (\ref{equation38'}).
Since the potential $(\alpha\sqrt{\xi})$ tends to $+\infty$ when
$\xi$ goes to $+\infty$, $E_1^\infty[\alpha\sqrt{\xi}]$ is reached by a positive function
$\psi\in S_\infty$:\[
-\psi''(\xi)+(\alpha \sqrt{\xi})\psi= E_1^\infty[\alpha \sqrt{\xi}]\psi.
\]
Setting $\xi=\alpha^\beta \zeta$ and
$\bar\psi(\zeta)=\sqrt{\alpha^\beta}\psi(\alpha^\beta \zeta)$ for an
arbitrary constant $\beta$, we
get\[
-\frac{1}{\alpha^{2\beta}}\bar\psi''(\zeta)+\alpha^{1+\frac{\beta}{2}}
\sqrt{\zeta}\bar\psi(\zeta)= E_1^\infty[\alpha
\sqrt{\xi}]\bar\psi(\zeta).
\]
By choosing $\beta$ such that
$-2\beta=1+\frac{\beta}{2}$, so that $\beta=\frac{-2}{5}$, we
obtain\[
-\bar\psi''(\zeta)+\sqrt{\zeta}\bar\psi(\zeta)= \alpha^{\frac{-4}{5}}E_1^\infty[\alpha \sqrt{\xi}]\bar\psi(\zeta),
\]
which implies that
$E_1^\infty[\alpha\sqrt{\xi}]=\alpha^{\frac{4}{5}}E_1^\infty[\sqrt{\xi}]$.
The proof of Proposition~\ref{proposition3.1} is achieved.

\subsection{Proof of Theorem \ref{mainth1}}

 In what follows, we will show that
(\ref{equation35}) admits a unique solution verifying (\ref{mod5}).
Indeed, the functional $J_0(.)$ is obviously continuous and strictly
convex on $\dot H^1_0(\mathbb{R}^+)$. To prove the existence of a
unique $U_0\in \dot H^1_0(\mathbb{R}^+)$ satisfying
\eqref{equation35}, it remains to verify that $J_0(.)$ is coercive
on $\dot H^1_0(\mathbb{R}^+)$. For this let $U\in \dot
H^1_0(\mathbb{R}^+)$ and write
$U(\xi)=\int_0^\xi U^\prime (t)dt$. This implies that $U(\xi)\leq
\|U^\prime\|_{L^2}\sqrt{\xi}$, and since $E_1^\infty[.]$ is an
increasing function one has $E_1^\infty[U]\leq
E_1^\infty[\|U^\prime\|_{L^2}\sqrt{\xi}]$. Applying
(\ref{equation38'}) with $\alpha=\|U^\prime\|_{L^2}$ we get
$E_1^\infty[U]\leq\|U^\prime\|_{L^2}^{\frac{4}{5}}E_1^\infty[\sqrt{\xi}]$, and finally we
have\[
J_0(U)\geq \frac{1}{2}\|U^\prime\|_{L^2}^2-E_1^\infty[\sqrt{\xi}].
\|U^\prime\|_{L^2}^{\frac{4}{5}}\xrightarrow[\|U\|_{H^1}\rightarrow +\infty]{}+\infty.
\]
Let us now prove that $U_0$ is a solution of the limit
problem (\ref{mod5}). Namely, we have to check that $E_1^\infty[U_0]$ is an eigenvalue.
To this aim, we first write the Euler--Lagrange equation
for~$U_0$:\begin{equation}
\label{equation44}
-U_0''=\frac{dE_1^\infty}{dU}[U_0]\geq 0.
\end{equation}
Therefore, $U_0$ is a concave function belonging to $\dot
H^1_0(\mathbb{R}^+)$. It is thus a continuous, increasing, and
positive function on $\mathbb{R}^+$, and
$\lim_{+\infty}U_0$ exists in $\overline{\mathbb{R}}$. It now remains to check that
 $E_1^\infty[U_0]<\lim_{+\infty}U_0$,  which will ensure that
$E_1^\infty[U_0]:=E_{1,0}$ is an eigenvalue with unique positive
eigenfunction $\psi_{1,0}\in S_{\infty}$ (see point two of
Proposition~\ref{proposition3.1}). We proceed by contradiction and
assume that $E_1^\infty[U_0]= \lim_{+\infty}U_0$.
Applying the third point of Proposition~\ref{proposition3.1}, one
obtains $\frac{dE_1^\infty}{dU}[U_0]=0$. In view of
(\ref{equation44}) and the fact that $U_0$ is a concave positive
function in $\dot H^1_0(\mathbb{R}^+)$, we deduce that $U_0=0$ and
$\min_{U\in \dot H^1_0(\mathbb{R}^+)}J_0(U)=0$. But a
simple rescaling argument shows that $J_0$ takes negative values, and
so its minimum is negative. To prove this claim, we fix a potential
$U$ in $\dot H^1_0(\mathbb{R}^+)$ such that $\int_0^{+\infty}|U^\prime|^2=1$,
$\lim_{+\infty}U=+\infty$ and let $\psi_1\in S_\infty$
be the eigenfunction corresponding to $E_1^\infty[U]$. For
$\varepsilon>0$,
 setting $U^\varepsilon(\xi)=\varepsilon^2U(\varepsilon \xi)$ and
$\psi_1^\varepsilon(\xi)=\sqrt{\varepsilon}\psi_1(\varepsilon \xi)$, we
have\[
-\frac{d^2\psi_1^\varepsilon}{d\xi^2}+U^\varepsilon \psi_1^\varepsilon(\xi)= \varepsilon^2E_1^\infty[U]\psi_1^\varepsilon(\xi),
\]
which implies that
$E_1^\infty[U^\varepsilon]=\varepsilon^2E_1^\infty[U]$. After
straightforward computations, we finally
obtain\begin{align*}
J_0(U^\varepsilon)= \frac{1}{2}\int_0^{+\infty}\left|\frac{dU^\varepsilon}{d\xi}\right|^2d\xi-E_1^\infty[U^\varepsilon]
&= \frac{\varepsilon^5}{2}\int_0^{+\infty}|U^\prime|^2d\xi-\varepsilon^2E_1^\infty[U] \\
&= -\varepsilon^2E_1^\infty[U]\left(1-\frac{\varepsilon^3}{2E_1^\infty[U]}\right),
\end{align*}
which is negative for $\varepsilon$ small enough. The proof of Theorem~\ref{mainth1} is
complete.

\section{Convergence analysis}\label{convergence}

The various models presented in the first section of this work are
all well posed. In this section, we shall estimate the difference
between their solutions in terms of $\eps$. Namely, we have to prove
estimates (\ref{equation18}) and (\ref{equation17}). The following
lemma will be useful.

\begin{lemma}\label{lemma4.2}
Let\/ $(U_0, E_{1,0}, \psi_{1,0})$ be the solution of the limit
problem\/ {\rm(\ref{mod5})}. There exist $a,b\in \mathbb{R}^+$ independent
of $\varepsilon$ such that for all $\varepsilon$ small we
have\begin{equation}
\label{equation45}
\|\psi_{1,0}\|_{L^2(M_\varepsilon,+\infty)}\leq a e^{-bM_\varepsilon}
\end{equation}
with $M_\varepsilon= \frac{1}{\varepsilon}$.\end{lemma}

\begin{proof}
We have $-\psi_{1,0}''+U_0\psi_{1,0}=E_{1,0}\psi_{1,0}$ such that
$\psi_{1,0}\geq 0$, $\psi_{1,0}(\xi)\rightarrow_{\xi\rightarrow +\infty} 0$ ($\psi_{1,0}\in H^1_0(\RR^+)$),
$E_{1,0}<\lim_{+\infty}U_0$, and $U_0$ increases to its limit at
$+\infty$. Then one can find two nonnegative constants $c$ and
$\delta $ independent of $\varepsilon$ such that for all
$\varepsilon$ small enough we
have\[
-\psi_{1,0}''(\xi)\leq -\delta \psi_{1,0}(\xi) \quad
\text{for } \xi\in [M_\varepsilon, +\infty[
\]
and $\psi_{1,0}(M_\varepsilon)\leq ce^{-\sqrt{\delta}M_\varepsilon}$. Let
$S(\xi)=ce^{-\sqrt{\delta}\xi}$ and $\psi=\psi_{1,0}-S$.  Then we
have $S''=\delta S$
and\begin{equation}
\label{varphi}
-\psi''(\xi)+\delta \psi(\xi)\leq 0, \qquad \xi\in[M_\varepsilon, +\infty[
\end{equation}
with\[
\psi(M_\varepsilon)\leq 0 \quad \text{and} \quad \psi(+\infty)= 0.
\]
By the maximum principle, one deduces that $\psi\leq 0$ on $[M_\eps,+\infty[$.
Thus\[
\psi_{1,0}(\xi)\leq ce^{-\sqrt{\delta}\xi}\quad
\text{on } [M_\varepsilon,+\infty[,
\]
which yields estimate (\ref{equation45}).\qquad\end{proof}

We begin by proving the second estimate \eqref{equation17} of
Theorem~\ref{mainth2}. For this we will compare (see Proposition~\ref{proposition4.1}) the potentials $U_0$ and $\tilde U_\eps$
solutions of \eqref{equation35} and \eqref{equation34}, respectively.
This will be done thanks to an idea consisting of the reformulation
of the problems \eqref{mod5} and \eqref{mod4} as minimization
problems whose unknown is the first eigenfunction. This is the
subject of the following remark.

\begin{remark}\label{remark4.1}
For $\phi\in S_{M_\varepsilon}$(see Definition\/~{\rm\ref{definition2.1}}),
where $M_\eps=\frac{1}{\eps}$, let us
set\begin{equation}
\label{equation47}
A_\varepsilon(\phi)=\int_0^{M_\varepsilon}|\phi^\prime(\xi)|^2d\xi +\frac{1}{2}\int_0^{M_\varepsilon}\!\int_0^{M_\varepsilon}\phi^2(\xi)\phi^2(\zeta)\min(\xi,\zeta)d\xi d\zeta
\end{equation}
and for
$\phi \in S_\infty$\begin{equation}
\label{equation46}
A_0(\phi)=\int_0^{+\infty}|\phi^\prime(\xi)|^2d\xi +\frac{1}{2}\int_0^{+\infty}\!\int_0^{+\infty}\phi^2(\xi)\phi^2(\zeta)\min(\xi,\zeta)d\xi d\zeta.
\end{equation}
The functional $A_\varepsilon$  satisfies $A_\eps(|\phi|)  = A_\eps(\phi)$ and the convexity
property\[
A_\eps\left(\sqrt{t\phi_1^2 + (1-t) \phi_2^2}\right) \leq t A_\eps(\phi_1) + (1-t) A_\eps(\phi_2)
\]
 for $t\in (0,1)$, the inequality being strict if\/ $|\phi_1|$
 and\/ $|\phi_2|$ are not proportional (these properties are also satisfied
by $A_0$). The functionals are obvious weakly lower semicontinuous
on their domain of definition, in such a way that the minimization
problems\begin{equation}
\label{equation49}
A_\varepsilon(\phi_\varepsilon)=\min_{\phi\in H^1_0(0,M_\varepsilon),\, \|\phi\|_{L^2}=1}A_\varepsilon(\phi)
\end{equation}
and\begin{equation}
\label{equation48}
A_0(\phi_0)=\min_{\phi\in H^1_0(\mathbb{R}^+),\, \|\phi\|_{L^2}=1}A_0(\phi)
\end{equation}
have unique positive solutions. The problems\/ {\rm(\ref{equation49})} and\/ {\rm(\ref{equation48})}
are equivalent, respectively, to\/ {\rm(\ref{mod4})} and\/ {\rm(\ref{mod5})}.
Indeed, the functions $\phi_\varepsilon$ and $\phi_0$
satisfy\begin{gather*}
-\phi_\varepsilon''+ U(\phi_\varepsilon)\phi_\varepsilon= \mu_\varepsilon \phi_\varepsilon \quad \text{on\/ }[0,M_\varepsilon],
\\ \noalign{\vspace*{-\jot}\vspace*{\abovedisplayskip}}
-\phi_0''+ U(\phi_0)\phi_0=\mu_0 \phi_0\quad \text{on\/ } \mathbb{R}^+,
\end{gather*}
where $\mu_\varepsilon$ (respectively, $\mu_0$) is the
Lagrange multiplicator
associated with the constraint\/ $\|\phi\|_{L^2}=1$ and $U(\phi_\varepsilon)$, $U(\phi_0)$ denote,
respectively,\[
U(\phi_\varepsilon)(\xi)=\int_0^{M_\varepsilon}|\phi_\varepsilon(\zeta)|^2\min(\xi,\zeta)d\zeta,\quad
U(\phi_0)(\xi)=\int_0^{+\infty}|\phi_0(\zeta)|^2\min(\xi,\zeta)d\zeta.
\]
In addition, since $\phi_\varepsilon$ and $\phi_0$ are
positive and the function $K(\xi,\zeta)=\min(\xi,\zeta)$ is the
kernel corresponding to the Laplacian in dimension one, we
have\[
 (U(\phi_\varepsilon),\mu_\varepsilon,\phi_\varepsilon)= (\tilde U_\varepsilon,\tilde E_{1,\varepsilon}, \tilde\psi_{1,\varepsilon})
 \quad \text{and\/}\quad(U(\phi_0),\mu_0,\phi_0)=(U_0,E_{1,0}, \psi_{1,0}).
\]\end{remark}\unskip

\begin{proposition}\label{proposition4.1}
The solutions $U_0$ and $\tilde U_\varepsilon$ of {\rm(\ref{equation35})}
and {\rm(\ref{equation34})}, respectively, verify
the following
estimate:\begin{equation}
\label{equation50}
\left\|\frac{d}{d\xi}(\tilde U_\varepsilon-U_0)\right\|_{L^2(0,M_\varepsilon)}^2=\mathcal{O}(e^{-\frac{c}{\eps}}),
\end{equation}
where $c$ is a strictly positive constant independent of
$\varepsilon$ and $M_\eps=\frac{1}{\eps}$. This yields estimate\/ {\rm(\ref{equation17})}.\end{proposition}

\begin{proof}
We start by comparing $A_0(\psi_{1,0})$ and
$A_\varepsilon(\tilde\psi_{1,\varepsilon})$. Let $\chi_\varepsilon
\in \mathcal{D}(0,+\infty)$ be such that $\chi_\varepsilon(\xi)=1$
on $[0,M_\varepsilon-1]$, $\chi_\varepsilon(\xi)=0$ on
$[M_\varepsilon, +\infty[$, and $0\leq\chi_\varepsilon\leq 1$. The
function $\chi_\varepsilon.\psi_{1,0}|_{(0,M_\varepsilon)}$ belongs
to $H^1_0(0,M_\varepsilon)$, and for $\varepsilon$ small we have
$\|\chi_\varepsilon.\psi_{1,0}\|_{L^2(0,M_\varepsilon)}\neq\nobreak 0$. Let
$\beta_\varepsilon=\|\chi_\varepsilon.\psi_{1,0}\|_{L^2(0,M_\varepsilon)}$.
Then we have $\frac{1}{\beta_\eps}\chi_\varepsilon.\psi_{1,0}\in S_{M_\varepsilon}$ and, with (\ref{equation45}),
$\beta_\varepsilon=1+O(e^{-\frac{c}{\eps}})$. Then, in view of
Remark~\ref{remark4.1}, the following inequalities can be straightforwardly
justified:\[
A_\varepsilon(\tilde\psi_{1,\varepsilon})\leq
A_\varepsilon\left(\frac{1}{\beta_\varepsilon}\chi_\varepsilon.\psi_{1,0}\right)
\leq A_\eps(\psi_{1,0})+ \mathcal{O}(e^{-\frac{c}{\eps}})
\leq A_0(\psi_{1,0})+ \mathcal{O}(e^{-\frac{c}{\eps}}),
\]
where $c$ is a  strictly positive constant independent of $\varepsilon$. Besides, we
have\[
A_0(\psi_{1,0})\leq A_0(\tilde\psi_{1,\varepsilon})=A_\varepsilon(\tilde\psi_{1,\varepsilon}).
\]
Here and in what follows, we still denote by
$\tilde\psi_{1,\varepsilon}$ the extension of
$\tilde\psi_{1,\varepsilon}$ by zero on $[M_\varepsilon,+\infty[$
when it is taken as a function on $\mathbb{R}^+$. Consequently, we
have\begin{equation}
\label{A_eps-A_0}
|A_\varepsilon(\tilde \psi_{1,\varepsilon})-A_0(\psi_{1,0})|=
A_\varepsilon(\tilde \psi_{1,\varepsilon})-A_0(\psi_{1,0})=
\mathcal{O}(e^{-\frac{c}{\eps}}).
\end{equation}
Furthermore, $A_0$ is uniformly convex on $S_\infty$ and
$\psi_{1,0}$ realizes its minimum. Then one can find a constant $c_0>0$ independent of $\varepsilon$ such
that\[
\|\psi_{1,0}-\tilde\psi_{1,\varepsilon}\|_{H^1(\mathbb{R}^+)}^2\leq c_0 |A_0(\psi_{1,0})-A_0(\tilde\psi_{1,\varepsilon})|.
\]
In addition, since $A_0(\tilde\psi_{1,\varepsilon})=A_\varepsilon(\tilde\psi_{1,\varepsilon})$ and with \eqref{A_eps-A_0}, one deduces
that\begin{equation}
\label{equation51}
\|\psi_{1,0}-\tilde\psi_{1,\varepsilon}\|_{H^1(\mathbb{R}^+)}^2=\mathcal{O}(e^{-\frac{c}{\eps}}).
\end{equation}
The potential $U_0-\tilde U_\varepsilon$
satisfies\[
-\frac{d^2}{d\xi^2}(U_0-\tilde U_\varepsilon)(\xi)= |\psi_{1,0}(\xi)|^2- |\tilde\psi_{1,\varepsilon}(\xi)|^2 \quad \text{on } [0,M_\varepsilon].
\]
Then, multiplying this equation by $U_0-\tilde U_\varepsilon$, one obtains after integration by
parts\begin{align*}
\left\|\frac{d}{d\xi}(U_0-\tilde U_\varepsilon)\right\|^2_{L^2(0,M_\varepsilon)}\leq {} & \frac{dU_0}{d\xi}(M_\eps)(U_0(M_\eps)-\tilde U_\varepsilon(M_\eps))\\
& {+}\> {\sup_{\xi\in[0,M_\eps]}(|U_0(\xi)-\tilde U_\varepsilon(\xi)|) \int_0^{M_\eps}\, (|\psi_{1,0}|^2-|\tilde\psi_{1,\eps}|^2)d\xi.}
\end{align*}
Moreover, in view of Remark~\ref{remark4.1}, we have, for every
$\xi\in [0,M_\eps]$,\begin{align*}
(U_0-\tilde U_\eps)(\xi)&= \int_0^{M_\eps}\,(|\psi_{1,0}|^2-|\tilde \psi_{1,\eps}|^2)(\zeta)\min(\xi,\zeta)d\zeta +\int_{M_\eps}^{+\infty}|\psi_{1,0}|^2\min(\xi,\zeta)d\zeta \\
&\leq M_\eps\left(\int_0^{M_\eps}\,(|\psi_{1,0}|^2-|\tilde \psi_{1,\eps}|^2)(\zeta)d\zeta +\int_{M_\eps}^{+\infty}|\psi_{1,0}|^2(\zeta)d\zeta \right)
\end{align*}
and $\frac{dU_0}{d\xi}(M_\eps)=-\int_{M_\eps}^{+\infty}\frac{d^2U_0}{d\xi^2}(\xi)d\xi=\int_{M_\eps}^{+\infty}|\psi_{1,0}|^2d\xi$.
Then\begin{align*}
\left\|\frac{d}{d\xi}(U_0-\tilde U_\eps)\right\|_{L^2(0,M_\eps)}^2 &\leq M_\eps\left(\int_0^{M_\eps}\,(|\psi_{1,0}|^2-|\tilde \psi_{1,\eps}|^2)(\zeta)d\zeta
+\int_{M_\eps}^{+\infty}|\psi_{1,0}|^2(\zeta)d\zeta \right)^2 \\
&\leq 2 M_\eps\left(\|\psi_{1,0}-\tilde \psi_{1,\eps}\|^2_{L^2(0,M_\eps)} +\|\psi_{1,0}\|_{L^2(M_\eps,+\infty)}^4\right),
\end{align*}
and with \eqref{equation51}\vadjust{\vspace*{1pt}}
and \eqref{equation45} one obtains \eqref{equation50}.
Moreover, with the change of variable $\tilde V_\eps(.)=\frac{1}{\eps^2}\tilde U_\eps(\frac{.}{\eps})$
one deduces
that\[
\left\|\frac{d}{d\xi}\left(\tilde V_\eps-\frac{1}{\eps^2}\tilde U_0\left(\frac{.}{\eps}\right)\right)\right\|_{L^2(0,1)}^2
=\frac{1}{\eps^5}\left\|\frac{d}{d\xi}(\tilde U_\eps-U_0)\right\|_{L^2(0,M_\eps)}^2=\mathcal{O}(e^{-\frac{c}{\eps}}),
\]
and then estimate \eqref{equation17} holds.\qquad\end{proof}

Let us now give the following result, which shows the
existence of a uniform gap between the first eigenvalue $\tilde
E_{1,\eps}:=E_1[\tilde U_\varepsilon]$ and the others $E_p[\tilde
U_\varepsilon]$.

\begin{lemma}\label{gap}
There exists a constant $G>0$, independent of $\varepsilon$, such
that\begin{equation}
\label{equation54}
E_p[\tilde U_\varepsilon]-\tilde E_{1,\eps}\geq G\quad \forall  p\geq 2.
\end{equation}\end{lemma}\unskip

\begin{proof}
Since $(E_p[\tilde U_\varepsilon])_{p\geq 1}$ is an increasing
sequence, it is sufficient to show (\ref{equation54}) only for
$p=2$. We argue by contradiction and suppose that $|E_2[\tilde
U_\varepsilon]-\tilde E_{1,\eps}|\rightarrow 0$ as $\varepsilon$
goes to zero. In view of Remark~\ref{remark4.1}, we have $\tilde
E_{1,\eps}=A_\eps(\tilde \psi_{1,\eps})$ and
$E_{1,0}=A_0(\psi_{1,0})$. Then, with \eqref{A_eps-A_0}, $|\tilde
E_{1,\eps}- E_{1,0}|=\mathcal{O}(e^{-\frac{c}{\eps}})$. We deduce
that $E_2[\tilde U_\varepsilon]$ and $\tilde E_{1,\eps}$ converge to
$E_{1,0}$ when $\eps \rightarrow 0$. Then the eigenfunctions
$\psi_2[\tilde U_\varepsilon]$ and $\tilde \psi_{1,\varepsilon}$,
corresponding, respectively, to $E_2[\tilde U_\varepsilon]$ and
$\tilde E_{1,\eps}$, prolonged by zero on
$[\frac{1}{\eps},+\infty)$, are bounded in $H^1_0(\mathbb{R}^+)$ with
respect to $\varepsilon$. There exists $\psi_1$ (respectively,
$\psi_2) \in H^1_0(\mathbb{R}^+)$ such that $\tilde
\psi_{1,\varepsilon}$ (respectively, $\psi_2[\tilde U_\varepsilon]$)
converges weakly in $H^1_{loc}(\mathbb{R}^+)$ to $\psi_1$
(respectively, $\psi_2$). By passing to the limit
$\varepsilon\rightarrow 0^+$ in $\mathcal{D}^\prime(0,+\infty)$, in the
equations\begin{gather*}
-\tilde \psi_{1,\varepsilon}''+\tilde U_\varepsilon \tilde \psi_{1,\varepsilon} =\tilde E_{1,\varepsilon} \tilde \psi_{1,\varepsilon},
\\ \noalign{\vspace*{-\jot}\vspace*{\abovedisplayskip}}
- (\psi_{2}[\tilde U_\varepsilon])''+\tilde U_\varepsilon \psi_{2}[\tilde U_\varepsilon] = E_{2}[\tilde U_\varepsilon] \psi_2[\tilde U_\varepsilon]
\end{gather*}
one deduces that $\psi_1$ and $\psi_2$ are two eigenfunctions
corresponding to $E_{1,0}$. In addition, we have
\begin{align*}
|J_{U_0}(\tilde \psi_{1,\varepsilon})-E_{1,0}|&= \left|\int_0^{+\infty}|\tilde \psi_{1,\eps}'|^2 +\int_0^{+\infty}U_0|\tilde \psi_{1,\eps}|^2-E_{1,0}\right|\\
&= \left|\int_0^{M_\eps}|\tilde \psi_{1,\eps}'|^2 +\int_0^{M_\eps}\tilde U_\eps|\tilde \psi_{1,\eps}|^2-E_{1,0} +\int_0^{+\infty}\,(U_0-\tilde U_\eps)|\tilde \psi_{1,\eps}|^2\right|\\
&\leq |\tilde E_{1,\eps}-E_{1,0}|+ \sup_{[0,M_\eps]}(|U_0-\tilde U_\eps|) \xrightarrow[\eps\rightarrow 0]{}0.
\end{align*}
Then $(\tilde \psi_{1,\varepsilon})$ (and similarly
$(\psi_2[\tilde U_\eps])$) is a minimizing sequence of ``$ E_{1,0}=\inf_{\psi\in S_\infty}J_{U_0}(\psi)$.'' Moreover, since
$E_{1,0}< \lim_{+\infty}U_0$ (see the proof of Theorem \ref{mainth1})
and applying Lemma \ref{lemmaA.1}, $(\tilde \psi_{1,\varepsilon})$
and $(\psi_2[\tilde U_\eps])$ (up to extraction of subsequences)
converge strongly in $L^2(\RR^+)$. Thus, since $\tilde
\psi_{1,\varepsilon}$ and $\psi_2[\tilde U_\eps]$ are two normalized
and orthogonal functions in $L^2(\RR^+)$ for any $\eps >0$, we
deduce that their limits when $\eps\rightarrow 0$, $\psi_1$ and
$\psi_2$, which are two eigenfunctions of $E_{1,0}$, are also normalized and orthogonal in $L^2(\RR^+)$. This
contradicts the fact that $E_{1,0}$ is a simple eigenvalue.\nobreak\qquad\end{proof}

\begin{proposition}\label{proposition4.2}
The potentials $U_\varepsilon$ and $\tilde U_\varepsilon$ solutions of\/ {\rm(\ref{mod3})}  and\/ {\rm(\ref{mod4})}
verify that\begin{equation}
\label{equation52}
\|U_\varepsilon-\tilde U_\varepsilon\|_{H^1(0,M_\varepsilon)}=\mathcal{O}(e^{-\frac{c}{\eps^2}}),
\end{equation}
where $c$ is a strictly positive constant independent of $\varepsilon$.
This gives estimate\/ {\rm(\ref{equation18})}.\end{proposition}

\begin{proof}
Recall first that $U_\varepsilon$ and $\tilde U_\varepsilon$ verify,
respectively, (\ref{equation33}) and (\ref{equation34}). To prove
estimate (\ref{equation52}), it is sufficient to compare the energies
$J_\varepsilon(U_\varepsilon)$ and $J_\varepsilon(\tilde U_\varepsilon)$ because we
have\begin{equation}
\label{equation53'}
\|U_\varepsilon-\tilde U_\varepsilon\|_{H^1(0,M_\varepsilon)}^2\leq c_0|J_\varepsilon(U_\varepsilon)-J_\varepsilon(\tilde U_\varepsilon)|,
\end{equation}
where $c_0$ is independent of $\varepsilon$. A
straightforward comparison gives the following
inequalities:\begin{equation}
\label{equation53}
\tilde J_\varepsilon(\tilde U_\varepsilon)\leq \tilde J_\varepsilon(U_\varepsilon)\leq J_\varepsilon(U_\varepsilon)
\leq J_\varepsilon(\tilde U_\varepsilon).
\end{equation}
Besides, we have $\tilde U_\varepsilon\geq 0$, and $E_p[.]$
is an increasing function; then $E_p[\tilde U_\varepsilon]\geq
E_p[0]=\eps^2p^2\pi^2 $. Moreover, since $\tilde E_{1,\eps}$
converges to $E_{1,0}$, which is then finite,
there exists a constant
$c_1>0$ independent of $\varepsilon$ such that
\begin{equation}\label{equation55}
E_p[\tilde U_\varepsilon]-\tilde E_{1,\eps}\geq \eps^2p^2\pi^2-c_1.
\end{equation}
Combining (\ref{equation54}) and (\ref{equation55}), one
finds $c_2>0$ and $c_3>0$ independent of $\varepsilon$ such
that\[
E_p[\tilde U_\varepsilon]-\tilde E_{1,\eps}\geq c_2\eps^2p^2\pi^2+ c_3\quad
\forall  p\geq 2.
\]
This implies
that\[
\sum_{p\geq 2}e^{-\frac{E_p[\tilde U_\varepsilon]-E_1[\tilde U_\varepsilon]}{\eps^2}}=\mathcal{O}(e^{-\frac{c_3}{\eps^2}}),
\]
and since\[
J_\varepsilon(\tilde U_\varepsilon)=\tilde J_\varepsilon(\tilde U_\varepsilon)
+\eps^2\log\,\Biggl(1 +\sum_{p\geq 2}  e^{-\frac{E_p[\tilde U_\varepsilon]-E_1[\tilde U_\varepsilon]}{\eps^2}}\Biggr)\,,
\]
we
obtain\[
J_\varepsilon(\tilde U_\varepsilon)=\tilde J_\varepsilon(\tilde U_\varepsilon)+ \mathcal{O}(\eps^2e^{-\frac{c_3}{\eps^2}}),
\]
which leads to  (\ref{equation52}) in view of (\ref{equation53'}) and (\ref{equation53}).\qquad\end{proof}

\section{Comments}\label{comments}

\subsection{Fermi--Dirac statistics}

It is more natural to consider
Fermi--Dirac statistics in the high density limit ($\eps\rightarrow
0$). Here we give some remarks and elements on the limit in this
case. The scaled occupation factor of the $p$th state with
Fermi--Dirac statistics is given
by\[
n_p^{FD}= f_{FD}(\mathcal{E}_{p}-\mathcal{E}_{F}),
\]
where $\mathcal{E}_{F}$ is the Fermi level and $f_{FD}$ is the Fermi--Dirac
distribution\begin{equation}
f_{FD}(u)= \log(1+e^{-u}),
\end{equation}
The scaled Boltzmann distribution function, however, is
given by $f_{B}(u)=e^{-u}$. The Poisson equation in model \eqref{mod1} can be written as
follows:\[
-\frac{d^2V}{d\xi^2}=\sum_{p=1}^{+\infty} f_{B(FD)}(\mathcal{E}_p-\mathcal{E}_F)|\varphi_p|^2
\]
under the following constraint on the Fermi
energy:\begin{equation}
\label{constr.}
\sum_{p=1}^{+\infty}f_{B(FD)}(\mathcal{E}_p-\mathcal{E}_F)=\frac{1}{\eps^3}.
\end{equation}
In the Boltzmann case, one can explicitly solve
\eqref{constr.} with respect to $\mathcal{E}_F$, and we
have\[
e^{\mathcal{E}_{F}}= \frac{1}{\eps^3 \sum_{p=1}^{+\infty}e^{-\mathcal{E}_{p}}},
\]
which yields \eqref{mod1}. The first remark we give in the
Fermi--Dirac case is that $e^{\mathcal{E}_{F}}$ cannot be
expressed explicitly in terms of $e^{-\mathcal{E}_{p}}$. The
analysis of the limit can, however, be extended to this case but with
technical complications that we have avoided in the Boltzmann
statistics case. When applying the change of variables
\eqref{equation1} and $\mathcal{E}_F=\frac{1}{\eps^2}\epsilon_F$
the intermediate problem \eqref{mod3} becomes in the Fermi--Dirac statistics
case\begin{equation}
\label{mod3FD}
{\left\{
\begin{array}{l@{}}
\ds -\frac{d^2\psi_{p}}{d\xi^2}+U\psi_{p}= E_{p}\psi_{p},\qquad \xi \in \left[0,\frac{1}{\varepsilon}\right], \\ \noalign{\vspace*{\jot}}
\ds \psi_{p}\in H^1\left(0,\frac{1}{\varepsilon}\right),\quad \psi_{p}(0)=0,\quad \psi_{p}\left(\frac{1}{\eps}\right)=0,\quad
\int_0^{\frac{1}{\varepsilon}} \psi_{p}\psi_{q}=\delta_{pq},\\ \noalign{\vspace*{\jot}}
\ds -\frac{d^2U}{d\xi^2}=\eps^3\sum_{p=1}^{+\infty} f_{FD}\left(\frac{E_p-\epsilon_F}{\eps^2}\right)|\psi_{p}|^2, \quad
\sum_{p=1}^{+\infty}f_{FD}\left(\frac{E_p-\epsilon_F}{\eps^2}\right)=\frac{1}{\eps^3}, \\ \noalign{\vspace*{\jot}}
\ds U(0)=0,\quad {dU\over d\xi}\left(\frac{1}{\varepsilon}\right)=0.
\end{array}
\right.}\hspace*{-\nulldelimiterspace}
\end{equation}
Since $f_{FD}$ is a regular, positive, and decreasing
function on $\RR$, the Schr\"odinger--Poisson system in a bounded
domain in the Fermi--Dirac case is well posed and can also be
expressed as an optimization problem; see the work of
Nier \cite{Nier1} in the unidimensional case and \cite{Nier3}~in higher
dimensions. More precisely, \eqref{mod3FD} is equivalent
to\begin{equation}
\label{minFD}
J_\eps(U_\eps)=\inf_{U\in H^{1,0}(0,\frac{1}{\eps})} J_\eps(U),
\end{equation}
where\begin{equation}
\label{functionalFD}
J_\varepsilon(U)=\frac{1}{2}\int_0^{\frac{1}{\eps}}|U^\prime|^2
- \varepsilon^3\sum_{p=1}^{+\infty}\,\Biggl[ f_{FD}\,\biggl(\frac{E_p[U]-\epsilon_F[U]}{\eps^2}\biggr)\,\epsilon_F[U]
-\eps^2\int_{\frac{E_p[U]-\epsilon_F[U]}{\eps^2}}^{+\infty}f_{FD}(u)du\Biggr]\,.
\end{equation}
Replacing $f_{FD}(.)$ by $f_B(u)=e^{-u}$, $J_\eps(.)$ is
nothing else but the functional \eqref{equation29} modulo a constant
independent of the variable $U$. The uniform gap showed in
Lemma~\ref{gap} remains correct. Then, as in the proof of
Proposition~\ref{proposition4.2}, there are two constants $c_1>0$ and $c_2>0$ such
that\[
(E_p-\epsilon_F)-(E_1-\epsilon_F)\geq c_1 \eps^2p^2\pi^2 + c_2\quad
\forall p\geq 2.
\]
Since $f_{FD}$ is a decreasing function and $\log(1+u)\sim u$ when $u\rightarrow 0^+$, then for all
$p\geq 2$\begin{align*}
\log\left(1+e^{-\frac{E_p-\epsilon_F}{\eps^2}}\right)&\leq \log\left(1+e^{-\frac{c_2}{\eps^2}}e^{-c_1p^2\pi^2}e^{-\frac{E_1-\epsilon_F}{\eps^2}}\right)
\leq c e^{-\frac{c_2}{\eps^2}}e^{-c_1p^2\pi^2}e^{-\frac{E_1-\epsilon_F}{\eps^2}}\\
&\leq c e^{-\frac{c_2}{\eps^2}}e^{-c_1p^2\pi^2} \log\left(1+e^{-\frac{E_1-\epsilon_F}{\eps^2}}\right),
\end{align*}
where $c >0$ is a general constant independent of $\eps$. This implies
that\[
\frac{\sum_{p\geq 2}\log\left(1+e^{-\frac{E_p-\epsilon_F}{\eps^2}}\right)}{\log\left(1+e^{-\frac{E_1-\epsilon_F}{\eps^2}}\right)}\leq
c \,\Biggl(\,\sum_{p\geq 2}e^{-c_1p^2\pi^2}\Biggr)\,e^{-\frac{c_2}{\eps^2}}.
\]
Thus, a formal analysis shows that, asymptotically when
$\eps\rightarrow 0$, \eqref{mod3FD} is close to a
Schr\"odinger--Poisson system with only the first energy level.
However, the rigorous analysis of the limit, $\eps\rightarrow 0$,
of \eqref{minFD}--\eqref{functionalFD} is more technically
complicated than the Boltzmann case for which the functional
$J_\eps$ has an explicit expression given by \eqref{equation29}.

\subsection{Boundary conditions and higher dimension}

The choice of
Neumann boundary condition at $z=1$ can be justified for modulation
doping devices (see~\cite{Bastard}) for which $z=1$ is in the bulk
of the semiconductor and the hypothesis of a vanishing electric
field is justified. This hypothesis also makes the analysis simple
because the boundary layer in the limit $\eps\rightarrow 0^+$ is
located at $z=0$. If $V$ satisfies Dirichlet boundary conditions,
then another boundary layer takes place at $z=1$. The analysis can
probably be extended to this case, but the first eigenvalue will have
asymptotically a multiplicity~2. The
multidimensional problem is more complicated, where the location of the electrons in the
boundary layer may depend on the geometry of the boundary. Such
problems have been noticed for the Schr\"odinger equation with a
magnetic field by \cite{Bonnaillie,Helffer} and are beyond the scope
of our work.\appendix

\section{Proof of Lemma \ref{lemmaA.1}}\label{appendix}
This appendix
is devoted to the proof of Lemma~\ref{lemmaA.1}. We will use the
concentration-compactness principle. This principle  is a general
method introduced by Lions~\cite{lions} to solve various
minimizing problems posed on unbounded domains. It is shown that all
minimizing sequences are relatively compact if and only if some
strict subadditivity inequalities hold. The proof is based upon a
lemma called the concentration-compactness lemma. For more details
on the principle, we refer the reader to~\cite{lions}. Let us begin by
recalling the concentration-compactness lemma.

\begin{lemma}[concentration-compactness lemma]\label{lemmaA2}
\begin{enumerate}
\item Let\/ $(\rho_n)_{n\geq 1}$ be a sequence in $L^1(\mathbb{R})$
satisfying $\rho_n\geq 0$ in\/ $\mathbb{R}$ and $\int_{\mathbb{R}}\rho_ndx=\lambda$ for a fixed $\lambda>0$. Then
there exists a subsequence\/ $(\rho_{n_k})_{k\geq 1}$ satisfying one
of the three following possibilities:
\begin{description}
\item{\rm(i)} (Compactness): there exists $y_k\in \mathbb{R}$ such
that\[
\forall \varepsilon>0,\quad
\exists R<+\infty,\quad
\int_{y_k+B_R} \rho_{n_k}(x)dx\geq \lambda-\varepsilon,
\]
where $B_R=\{x\in \mathbb{R};${\rm~}$|x|\leq R\}$.
\item{\rm(ii)} (Vanishing):
\[
\lim_{k\rightarrow +\infty}\sup_{y\in \mathbb{R}}\int_{y+B_R}\rho_{n_k}(x)dx =0.
\]
\item{\rm(iii)} (Dichotomy):
There exists $\alpha \in \mathopen{]}0,\lambda[$ such that for all $\varepsilon
>0$, there exist $k_0\geq 1$ and $\rho_k^1, \rho_k^2 \in L^1_+(\mathbb{R})$ such that for
$k\geq k_0$\[
\|\rho_{n_k}-(\rho_k^1+\rho_k^2)\|_{L^1}\leq \varepsilon,\quad
\biggl|\int_{\mathbb{R}}\rho_k^1 dx - \alpha\biggr|\leq \varepsilon,\quad
\biggl|\int_{\mathbb{R}}\rho_k^2 dx -(\lambda- \alpha)\biggr|\leq \varepsilon,
\]
where $\rho_k^1$ has compact support and\/ $\operatorname{dist}(\operatorname{supp} (\rho_k^1),
\operatorname{supp}(\rho_k^2)) \rightarrow_k +\infty$.
\end{description}
\item If $\rho_n= |u_n|^2$ with $u_n$ bounded in $H^1(\mathbb{R})$,
there exists a subsequence\/ $(\rho_{n_k})$ such that either
compactness\/~{\rm(i)}, vanishing\/~{\rm(ii)}, or dichotomy\/~{\rm(iii)} occurs
as follows: there exists $\alpha \in \mathopen{]}0,\lambda[$
such that for all $\varepsilon >0$ there exist $k_0\geq 1$, $u_k^1, u_k^2$ bounded in $H^1(\mathbb{R})$ satisfying for
$k\geq k_0$\[
{\left\{
\begin{array}{l@{}}
\ds \|u_{n_k}-(u_k^1+u_k^2)\|_{L^2}\leq\delta(\varepsilon) \xrightarrow[\varepsilon\rightarrow 0]{}0, \\  \noalign{\vspace*{\jot}}
\ds \left| \int_{\mathbb{R}}|u_k^1|^2dx-\alpha\right|\leq\varepsilon,\qquad \left|\int_{\mathbb{R}}|u_k^2|^2dx-(\lambda-\alpha)\right|\leq\varepsilon , \\ \noalign{\vspace*{\jot}}
\ds \operatorname{dist}\left(\operatorname{supp}(u_k^1), \operatorname{supp}(u_k^2)\right)\xrightarrow[k\rightarrow+\infty]{}+\infty, \\ \noalign{\vspace*{\jot}}
\ds   \liminf_{k} \int_{\mathbb{R}}\, \{|\nabla u_{n_k}|^2-  |\nabla u_k^1|^2-|\nabla u_k^2|^2 \}dx\geq 0 .
\end{array}
\right.}\hspace*{-\nulldelimiterspace}
\]
\end{enumerate}\end{lemma}

First, we need to give some notation. For $V\in
\dot H^1_0(\mathbb{R}^+)$ and $\varepsilon \in \mathbb{R}^+$, we
define\begin{equation}
\label{equationA1}
I_\varepsilon =\inf\left\{J_V(\varphi),\quad \varphi\in H^1_0(\mathbb{R}^+),\quad \int_0^{+\infty}\varphi^2=\varepsilon\right\},
\end{equation}
where  $J_V(\varphi)= \int_0^{+\infty}|\varphi^\prime|^2 +\int_0^{+\infty}V\varphi^2$,
and\begin{equation}
\label{equationA2}
I_\varepsilon^\infty=\inf\left\{J_\infty(\varphi), \quad \varphi\in H^1_0(\mathbb{R}^+),\quad \int_0^{+\infty}\varphi^2=\varepsilon\right\},
\end{equation}
where $J_\infty(\varphi)= \int_0^{+\infty}|\varphi^\prime|^2 +V^\infty\int_0^{+\infty}\varphi^2$ and
$V^\infty=\liminf_{+\infty}V$.

\begin{lemma}
Let $V\in \dot H^1_0(\mathbb{R}^+)$ such that
$E_1^\infty[V]<\liminf_{+\infty}V$. Then the following strict subadditivity inequality
holds:\begin{equation}
\label{equationA3}
I_\varepsilon<I_\alpha+I_{\varepsilon-\alpha}^\infty\quad \forall\,
0<\alpha<\varepsilon.
\end{equation}\end{lemma}\unskip

\begin{proof}
Take $\varphi\in H^1_0(\mathbb{R}^+)$ such that
$\int_0^{+\infty}\varphi^2=\varepsilon$ and let $\psi=\frac{\sqrt{\alpha}}{\sqrt{\varepsilon}}\varphi$. Then $\psi\in
H^1_0(\mathbb{R}^+)$, $\int_0^{+\infty}\psi^2=\alpha$,
and $J_V(\psi)=\frac{\alpha}{\varepsilon}J_V(\varphi)$. This
implies that $\varepsilon I_\alpha\leq \alpha I_\varepsilon$ for any
arbitrary $\alpha >0$ and $\eps>0$. We also  deduce that $\varepsilon
I_\alpha = \alpha I_\varepsilon$ (and similarly $\varepsilon
I_\alpha^\infty = \alpha I_\varepsilon^\infty$) for any $\eps,\alpha
>0$. In particular, if $\alpha =1$, $\eps I_1=I_\eps$ (or $\eps
I_1^\infty=I_\eps^\infty$) for any $\eps >0$. Moreover, by
definition  we have $I_1=E_1^\infty[V]<V^\infty$. Then, for all $\varphi\in
H^1_0(\mathbb{R}^+)$ such that
$\int_0^{+\infty}\varphi^2=1$, we have $I_1<\int_0^{+\infty}|\varphi^\prime|^2+V^\infty$. This implies that
$I_1<I_1^\infty$, and by multiplying by $(\varepsilon-\alpha)$, which
is positive if $0<\alpha<\varepsilon$, one obtains
$\eps I_1-\alpha I_1< (\eps-\alpha) I_1^\infty$ and inequality
(\ref{equationA3}) holds.\qquad\end{proof}

{\it Proof of Lemma\/ {\rm\ref{lemmaA.1}}}.
Applying the
concentration-compactness lemma for  $(\rho_n)_n$:
$\rho_n(x)=|\psi_n(x)|^2$ on $\mathbb{R}^+$ and zero elsewhere,
there exists a subsequence $(\rho_{n_k})_k$ satisfying one of the
three cases given by Lemma~\ref{lemmaA2}. If vanishing~(ii) occurs, i.e.,
if\[
\lim_{k\rightarrow +\infty}\sup_{y\in\mathbb{R}} \int_{y+B_R}\rho_{n_k}(x)dx =0\quad \forall R\geq 0,
\]
which implies
that\[
\lim_{k\rightarrow +\infty}\sup_{y\in\mathbb{R}^+} \int_y^{y+R}|\psi_{n_k}(x)|^2dx=0\quad \forall R\geq 0,
\]
then, for all $\varepsilon>0$ small enough, one can find a
sequence $(R_k)_k$ of increasing positive real such that
$\int_0^{R_k}|\psi_{n_k}(x)|^2dx\leq \varepsilon$ for
all $k$. We
have\begin{eqnarray*}
\int_0^{+\infty}|\psi_{n_k}^\prime|^2dx+\int_0^{+\infty}V\psi_{n_k}^2dx&=& \int_0^{+\infty}|\psi_{n_k}^\prime|^2dx+\int_0^{R_k}V\psi_{n_k}^2dx+\int_{R_k}^{+\infty}V\psi_{n_k}^2dx\\
&\geq&\int_0^{+\infty}|\psi_{n_k}^\prime|^2dx-\|V\|_\infty.\varepsilon+(V^\infty-\varepsilon) \int_{R_k}^{+\infty}\psi_{n_k}^2dx.
\end{eqnarray*}
This implies that there exists $\delta(\varepsilon)$,
tending to zero when $\varepsilon\rightarrow 0$, such
that\[
J_V(\psi_{n_k})\geq J_\infty(\psi_{n_k})-\delta(\varepsilon)\geq I_1^\infty-\delta(\varepsilon).
\]
Now let $k$ go to $+\infty$ and $\varepsilon$ to zero.  Then we
obtain\[
I_1\geq I_1^\infty,
\]
which contradicts the strict subadditivity inequality (\ref{equationA3}).

Now we assume that $(\rho_{n_k})_k$ verifies the dichotomy
case; i.e., there exists $\alpha \in \mathopen{]}0,1[$ such that for all
$\varepsilon >0$ there exist $k_0\geq 1$, $\psi_k^1, \psi_k^2$
bounded in $H^1(\mathbb{R}^+)$ satisfying for
$k\geq k_0$\[
{\left\{
\begin{array}{l@{}}
\ds  \|\psi_{n_k}-(\psi_k^1+\psi_k^2)\|_{L^2}\leq\delta(\varepsilon)  \xrightarrow[\varepsilon\rightarrow 0]{}0, \\ \noalign{\vspace*{\jot}}
\ds \left| \int_0^{+\infty}|\psi_k^1(x)|^2dx-\alpha\right|\leq\varepsilon,\qquad \left|\int_0^{+\infty}|\psi_k^2(x)|^2dx-(1-\alpha)\right|\leq\varepsilon ,\\ \noalign{\vspace*{\jot}}
\ds \operatorname{dist}\left(\operatorname{supp}(\psi_k^1) ,   \operatorname{supp}(\psi_k^2)\right)\xrightarrow[k\rightarrow+\infty]{}+\infty,\\ \noalign{\vspace*{\jot}}
\ds   \liminf_{k} \underbrace{\int_{\mathbb{R}^+}\,\{|\nabla \psi_{n_k}|^2  -|\nabla \psi_k^1|^2-|\nabla \psi_k^2|^2 \}dx}_{\gamma_k}\geq 0 .
\end{array}
\right.}\hspace*{-\nulldelimiterspace}
\]
One can write
(see~\cite{lions})\[
\psi_{n_k}=\psi_k^1+\psi_k^2+\varphi_k,\quad \text{where}\,\quad \psi_k^1\psi_k^2= \psi_k^1\varphi_k=\psi_k^2\varphi_k=0\,\quad \text{a.e.},
\]
and without loss of generality, we suppose that
$\operatorname{supp}(\psi_k^2)\subset [R_k,+\infty[$, where $R_k$ tends to $+\infty$
with $k$. This implies
that\begin{align*}
\int_0^{+\infty}V\psi_{n_k}^2&= \int_0^{+\infty}V|\psi_k^1|^2+ \int_0^{+\infty}V|\psi_k^2|^2+\int_0^{+\infty}V|\varphi_k|^2 \\
&\geq \int_0^{+\infty}V|\psi_k^1|^2+(V^\infty-\varepsilon) \int_0^{+\infty}|\psi_k^2|^2-\|V\|_\infty.\delta(\varepsilon)
\end{align*}
and\[
\int_0^{+\infty}|\psi_{n_k}^\prime|^2+\int_0^{+\infty} V\psi_{n_k}^2\geq
\int_0^{+\infty}|\psi_{n_k}^\prime|^2+\int_0^{+\infty} V|\psi_k^1|^2+V^\infty\int_0^{+\infty}|\psi_k^2|^2-\delta(\varepsilon).
\]
Hence,\begin{equation}
\label{equationA4}
J_V(\psi_{n_k})\geq\gamma_k +J_V(\psi_k^1)+J_\infty(\psi_k^2)-\delta(\varepsilon).
\end{equation}
Besides, let $\alpha_k= \int_0^{+\infty}|\psi_k^1(x)|^2dx$,
$\beta_k=\int_0^{+\infty}|\psi_k^2(x)|^2dx$. For all fixed $\varepsilon>0$,
the sequences $(\alpha_k)_k$ and $(\beta_k)_k$ are bounded in
$\mathbb{R}^+$. There are subsequences, still denoted by
$(\alpha_k)_k$ and $(\beta_k)_k$, which converge in $\mathbb{R}^+$
to $\alpha_\eps$ and $\beta_\varepsilon$, respectively, where
$\alpha_\eps$ and $\beta_\varepsilon$ belong to $\mathbb{R}^+$ such
that\begin{equation}
\label{equationA5}
|\alpha_\eps-\alpha|\leq \varepsilon\quad \text{and}\quad
|\beta_\varepsilon-(1-\alpha)|\leq \varepsilon.
\end{equation}
Inequality (\ref{equationA4})
yields\[
J_V(\psi_{n_k})\geq \gamma_k + I_{\alpha_k}+I_{\beta_k}^\infty-\delta(\varepsilon).
\]
Taking the $\liminf_k$ of the last inequality and
letting $\varepsilon$ tend to zero, we obtain in view of
(\ref{equationA5}) and the fact that
$\liminf_k \gamma_k\geq 0$\[
I_1\geq I_\alpha+I_{1-\alpha}^\infty,
\]
which contradicts the strict subadditivity inequality (\ref{equationA3}).

Consequently, the sequence $(\rho_{n_k})_k$ verifies the
compactness case of the con\-cen\-tra\-tion-compactness 
lemma which yields
straightforwardly that the minimizing sequence $(\psi_n)_n$ is
relatively compact in $L^2(\mathbb{R}^+)$.

\subsection*{Acknowledgment}

The authors would like to thank Claude Le~Bris for fruitful discussions and in
particular for having pointed out the variational formulation of
Remark~\ref{remark4.1}.



\begin{thebibliography}{99}

\bibitem{AFS}
{\sc T. Ando, A. B. Fowler, and F. Stern},
{\it Electronic properties of two-dimensional systems},
Rev. Modern Phys., 54 (1982), pp. 437--672.

\bibitem{Bastard}
{\sc G. Bastard},
{\it Wave Mechanics Applied to Semiconductor Heterostructures},
Les \smash{\'Editions} de Physiques, Les Ulis Cedex, France, 1996.

\bibitem{nba}
{\sc N. Ben Abdallah},
{\it On a multidimensional Schr\"odinger-Poisson scattering model for semiconductors},
J. Math. Phys., 41 (2000), pp. 4241--4261.

\bibitem{nba-degmar}
{\sc N. Ben Abdallah, P. Degond, and P. Markowich},
{\it On a one-dimensional Schr\"odinger-Poisson scattering model},
Z. Angew. Math. Phys., 48 (1997), pp. 135--155.

\bibitem{BM}
{\sc N. Ben Abdallah and F. M\'ehats},
{\it On a Vlasov-Schr\"odinger-Poisson model},
Comm. Partial Differential Equations, 29 (2004), pp. 173--206.

\bibitem{BMO}
{\sc N. Ben Abdallah, F. M\'ehats, and O. Pinaud},
{\it Adiabatic approximation of the Schr\"odinger--Poisson system with a partial confinement},
SIAM J. Math. Anal., 36 (2005), pp. 986--1013.

\bibitem{Bonnaillie}
{\sc V. Bonnaillie},
{\it On the fundamental state energy for a Schr\"odinger operator with magnetic field in domains with corners},
Asymptot. Anal., 41 (2005), pp. 215--258.

\bibitem{Castella}
{\sc F. Castella},
{\it $L^2$ solutions to the Schr\"odinger-Poisson system: Existence, uniqueness, time behaviour, and smoothing effects},
Math. Models Methods Appl. Sci., 7 (1997), pp. 1051--1083.

\bibitem{CBL}
{\sc I. Catto, C. Le Bris, and P.-L. Lions},
{\it Mathematical Theory of Thermodynamic Limits: Thomas-Fermi Type Models},
Clarendon Press, Oxford, UK, 1998.

\bibitem{CAZ}
{\sc T. Cazenave},
{\it Semilinear Schr\"odinger Equations}, Courant Lect. Notes Math. 10,
AMS, Providence, RI, 2003.

\bibitem{HFB}
{\sc J. Harris, C. T. Foxon, K. Barnham, D. Lacklison, J. Hewett, and C. White},
{\it Two-dimensional electron gas structures with mobilities in excess of\/ $3.10^6 cm^2 V^{-1}S^{-1}$},
J. Appl. Phys., 61 (1987), pp. 1219--1221.

\bibitem{Helffer}
{\sc B. Helffer and A. Mohamed},
{\it Semiclassical analysis for the ground state energy of a Schr\"o\-dinger operator with magnetic wells},
J. Funct. Anal., 138 (1996), pp. 40--81.

\bibitem{IZL}
{\sc R. Illner, P. F. Zweifel, and H. Lange},
{\it Global existence, uniqueness and asymptotic behavior of solutions of the Wigner-Poisson and Schr\"odinger-Poisson systems},
Math. Methods Appl. Sci., 17 (1994), pp. 349--376.

\bibitem{KaiserRehberg}
{\sc H.-C. Kaiser and J. Rehberg},
{\it On stationary Schr\"odinger-Poisson equations modelling an electron gas with reduced dimension},
Math. Methods Appl. Sci., 20 (1997), pp. 1283--1312.

\bibitem{weierstrass}
{\sc H.-C. Kaiser and J. Rehberg},
{\it About a stationary Schr\"odinger-Poisson system with Kohn-Sham potential in a bounded two- or three-dimensional domain},
Nonlinear Anal., 41 (2000), pp. 33--72.

\bibitem{Kato}
{\sc T. Kato},
{\it Perturbation Theory for Linear Operators},
Springer-Verlag, Berlin, Heidelberg, 1966.

\bibitem{lions}
{\sc P. L. Lions},
{\it The concentration-compactness principle in the calculus of variations. The locally compact case, part\/ {\rm 1}},
Ann. Inst. H. Poincar\'e Anal. Non Lin\'eaire, 1 (1984), pp. 109--145.

\bibitem{Mark}
{\sc P. A. Markowich},
{\it Boltzmann distributed quantum steady states and their classical limit},
Forum Math., 6 (1994), pp. 1--33.

\bibitem{MRS}
{\sc P. A. Markowich, C. A. Ringhofer, and C. Schmeiser},
{\it Semiconductor Equations},
Springer-Verlag, Vienna, 1990.

\bibitem{Nier1}
{\sc F. Nier},
{\it A stationary Schr\"odinger-Poisson system arising from the modelling of electronic devices},
Forum Math., 2 (1990), pp. 489--510.

\bibitem{Nier2}
{\sc F. Nier},
{\it Schr\"odinger-Poisson systems in dimension $d \leq 3$: The whole-space case},
Proc. Roy. Soc. Edinburgh Sect. A, 123 (1993), pp. 1179--1201.

\bibitem{Nier3}
{\sc F. Nier},
{\it A variational formulation of Schr\"odinger-Poisson systems in dimension $d \leq 3$},
Comm. Partial Differential Equations, 18 (1993), pp. 1125--1147.

\bibitem{Nier4}
{\sc F. Nier},
{\it The dynamics of some quantum open systems with short-range nonlinearities},
Nonlinearity, 11 (1998), pp. 1127--1172.


\bibitem{Pin}
{\sc O. Pinaud},
{\it Adiabatic approximation of the Schr\"odinger-Poisson system with a partial confinement: The stationary case},
J. Math. Phys., 45 (2004), pp. 2029--2050.

\bibitem{PN}
{\sc E. Polizzi and N. Ben Abdallah},
{\it Self-consistent three dimensional model for quantum ballistic transport in open systems},
Phys. Rev. B, 66 (2002), pp. 245301--245309.

\bibitem{nbapolizzi}
{\sc E. Polizzi and N. Ben Abdallah},
{\it Subband decomposition approach for the simulation of quantum electron transport in nanostructures},
J. Comput. Phys., 202 (2005), pp. 150--180.

\bibitem{P}
{\sc E. Polizzi},
{\it Mod\'elisation et simulation num\'eriques du transport quantique balistique dans les nanostructures semi-conductrices},
Ph.D. thesis, INSA, Toulouse, France, 2001.

\bibitem{Poschel}
{\sc J. Poschel and E. Trubowitz},
{\it Inverse Spectral Theory},
Academic Press, Boston, 1987.

\bibitem{Reed Simon}
{\sc M. Reed and B. Simon},
{\it Methods of Modern Mathematical Physics},
Academic Press, New York, San Francisco, London, 1975.

\bibitem{V}
{\sc B. Vinter},
{\it Subbands and charge control in a two-dimensional electron gas field effect transistor},
Appl. Phys. Lett., 44 (1987), pp. 307--309.

\end{thebibliography}
\end{document}